\RequirePackage{algorithm}
\RequirePackage{algpseudocode}

\documentclass[hidelinks,onefignum,onetabnum]{siamart250211}

\usepackage[T1]{fontenc}
\usepackage[american]{babel}
\usepackage{csquotes}
\usepackage{amsmath}
\usepackage{amsfonts}
\usepackage{amssymb}
\usepackage{mathtools}
\usepackage{xcolor}
\usepackage{algorithm}
\usepackage{algpseudocode}
\usepackage{nicematrix}
\usepackage{placeins}
\usepackage[export]{adjustbox}
\usepackage{multirow}
\usepackage{subcaption}
\usepackage{tikz}
\usetikzlibrary{decorations.pathreplacing,calligraphy}

\tikzset{
brace/.style = {decorate,
            decoration={calligraphic brace, amplitude=2pt,
            pre =moveto, pre  length=1pt,
            post=moveto, post length=1pt,
            raise=#1},
            thick,
            pen colour=black},
brace/.default=0pt
}

\crefname{line}{line}{lines}

\usepackage[left=3cm, right=3cm]{geometry}

\usepackage{xargs}
\usepackage[]{todonotes}
\newcommandx{\fixme}[2][1=]{\todo[linecolor=red,backgroundcolor=red!25,bordercolor=red,#1]{#2}}
\newcommandx{\change}[2][1=]{\todo[linecolor=blue,backgroundcolor=blue!25,bordercolor=blue,#1]{#2}}
\newcommandx{\info}[2][1=]{\todo[linecolor=green,backgroundcolor=green!25,bordercolor=green,#1]{#2}}
\newcommandx{\improvement}[2][1=]{\todo[linecolor=orange,backgroundcolor=orange!25,bordercolor=orange,#1]{#2}}
\newcommandx{\thiswillnotshow}[2][1=]{\todo[disable,#1]{#2}}

\usepackage{hyperref}
\hypersetup{
    colorlinks,
    citecolor=black,
    filecolor=black,
    linkcolor=black,
    urlcolor=black,
    hidelinks
}

\ifpdf
\hypersetup{
  pdftitle={A Parallel-in-Time Combination Method for Parabolic Problems},
  pdfauthor={M. Griebel, M. A. Schweitzer and L. Troska}
}
\fi

\headers{A Parallel-in-Time Combination Method for Parabolic Problems}{M. Griebel, M. A. Schweitzer and L. Troska}

\title{A Parallel-in-Time Combination Method for Parabolic Problems} 
\author{Michael Griebel\thanks{Institut f\"ur Numerische Simulation, Universit\"at Bonn,
 Friedrich-Hirzebruch-Allee 7 53115 Bonn and Fraunhofer Institut f\"ur Algorithmen und Wissenschaftliches Rechnen SCAI, Schloss Birlinghoven, 53757 Sankt Augustin} \and Marc Alexander Schweitzer\footnotemark[1] \and Lukas Troska\footnotemark[1]}

\date{}

\crefname{paragraph}{Section}{Sections}

\renewcommand{\vec}[1]{\mathbf{#1}}

\DeclarePairedDelimiter{\norm}{\lVert}{\rVert}

\DeclareMathOperator{\diag}{diag}
\DeclareMathOperator{\blockdiag}{blockdiag}

\newenvironment{remark}{\par\medskip\noindent Remark:}{\par
\medskip}

\begin{document}

\maketitle
\begin{abstract}
	In this article, we present a parallel discretization and solution method for parabolic problems with a higher number of space dimensions. It consists of a parallel-in-time approach \cite{PinT2024} using the multigrid reduction-in-time
	algorithm \texttt{MGRIT} \cite{Falgout2014} with its implementation in the library \texttt{XBraid} \cite{xbraid-package}, the sparse grid combination method \cite{Bungartz.Griebel:2004,Griebel1992} for discretizing the resulting elliptic problems in space, and a domain decomposition method \cite{Griebel2023} for each of the subproblems in the combination method based on the space-filling curve approach. As a result, we obtain a highly parallel solver with good speedup and scale-up qualities, which is well-suited for parabolic problems with up to six space dimensions. 

	We describe our new parallel approach and show its superior parallelization properties for the heat equation, the chemical master equation and some exemplary stochastic differential equations.
\end{abstract}



\section{Introduction}
We consider the numerical approximation of parabolic partial differential equations
\begin{equation}\label{pdgl}
	\frac{\partial u}{\partial t} + {\cal L}u = f \quad \text{in } \Omega\times(T_\text{start}, T_\text{end}],
\end{equation}
with $\Omega \subset \mathbb{R}^d$ and linear elliptic differential operator $\cal L$, where we are especially interested in the case of larger spatial dimensions $d$. Such problems arise, for example, from stochastic differential equations, where modeling with the probability density function of the underlying Markov process leads to the well-known general Kolmogorov forward or the Fokker-Planck equation. It is used, for example, to describe the particle velocities in gases and liquids in statistical mechanics and thermodynamics, the reactions of concentrations of reactants in a chemical system,  the dynamics of friction and wear of mechanical systems in engineering, the evolution of populations of species in cell biology, genetics  and ecology, the firing rate of neurons in neuroscience, the dynamics of quantum mechanical systems through the Wigner function, the evolution of asset prices or volatility for option pricing and portfolio optimization in financial mathematics or the modeling of cognitive processes in behavioral science and decision making in psychology, to name just a few applications.

For its numerical approximation a common approach is to employ the method of lines, where the problem is first discretized in space and the resulting system of ordinary differential equations is then treated by a  numerical integration approach such as the Runge-Kutta method. An alternative and more common approach is the Rothe method, where the problem is first discretized in time, resulting in a sequence of elliptic problems, each of which is then discretized in space by a finite element or finite difference method. The latter approach uses a sequential time stepping scheme, where within each time step an elliptic subproblem must be solved, where the operator and right hand side depend on $\cal  L$, $f$, and the respective time stepping method.

There are two major problems here: First, when it comes to parallelization, in Rothe's approach the use of large-scale parallel systems is usually limited to parallel elliptic solvers, i.e.\ only the spatial dimensions are considered for parallelization  and the sequence of time steps is processed only sequentially, which typically requires very long runtimes to achieve the necessary accuracy. Note that for sequential time integrators in the method of lines, one typically also only considers spatial parallelization and thus encounters similar runtimes. To address this problem, we resort to a parallel-in-time approach \cite{PinT2024}, namely the multigrid reduction-in-time algorithm \texttt{MGRIT} \cite{Falgout2014} via its implementation in the library \texttt{XBraid} \cite{xbraid-package}.
It gives us a parallelization of the whole problem for the time coordinate, which can be considered as a {\em first scale of parallelization}. 
Second, for larger spatial dimensions $d$ we encounter the well-known curse of dimension, either in the size of the ODE system for the method of lines or in the system of equations that arise from the discretization of the spatial problems in Rothe's method. There, a space discretization on a uniform grid with mesh size $2^L$ in each coordinate direction is usually used, with associated level $L$. 
Thus, the cost scales as $O(2^{Ld})$, i.e.\ it basically scales exponentially with $d$, which makes any conventional discretization approach impossible for $d>3$, both for the method of lines and for Rothe's approach.

\begin{figure}
    \centering
    \includegraphics[valign=m, scale=0.7]{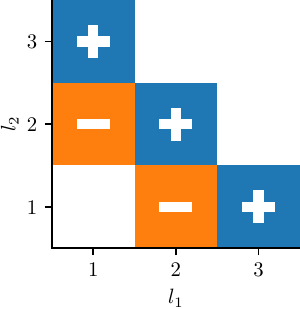}
    \includegraphics[valign=m]{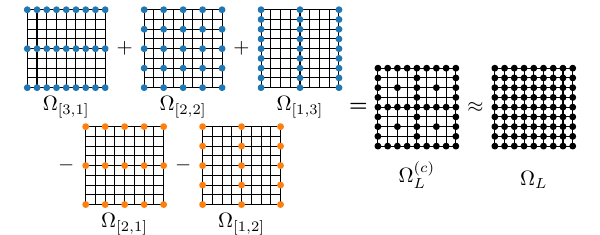}
    \caption{The combination method in two dimensions for level $L=3$.}\label{fig:combination_technique}
\end{figure}

In this article we will focus on Rothe's method, but instead of a uniform discretization in space, we will use a sparse grid approach \cite{Bungartz.Griebel:2004} to discretize the elliptic problems at each time step in the form of the so-called combination method \cite{Griebel1992}. This alleviates the curse of dimension and practically allows values of $d$ up to six. It exploits the fact that the resulting elliptic problems in the time step sequence (right hand side and boundary conditions permitting) generally have a smooth solution in each time step (perhaps except for the first one). This is due to the smoothing effect of $\cal L$ over time. The sparse grid combination  method involves solving the elliptic problems discretized as several different subproblems on (mostly) anisotropic grids $\Omega_l$, with mesh sizes $h_l=(2^{-l_1}, \ldots, 2^{-l_d})$, where 
\begin{equation}
l=(l_1,\ldots,l_d) \quad \mbox{ with } l_1+\ldots l_d =L+d-1-w, \quad w=0,\ldots, d-1,
\end{equation}
and $L$ is the level of the corresponding classical discretization with uniform mesh size $(2^{-L},\ldots,2^{-L})$ in each time step of Rothe's method, which would be prohibitive due to the curse of dimension. These solutions of the subproblems on these coarse anisotropic meshes are then combined to obtain a solution on the sparse grid associated to level $L$. We refer to \cref{fig:combination_technique} for a two dimensional example of the different subproblems and their combination for $L = 3$. For practical reasons, we also consider a coarsest minimal level $l>L_0:=(L_{0,1}, \ldots, L_{0,d})$. All these different subproblems of the combination method per each time step can be solved completely in parallel, which is a {\em second scale of parallelization}. Only their solutions have to be combined properly, which involves communication and some sequential computations.
Third, for the solution of each of the subproblems on the levels $l$ of the combination method and for each time step, we use a domain decomposition method, which is based on the space filling curve approach. It was already developed and analyzed in detail in \cite{Griebel2023}. Thus, the resulting subdomain problems per subproblem (plus a coarse scale problem) of the combination method for each time step can also be treated in parallel, which involves a {\em third scale of parallelization}. 

\begin{figure}
    \centering
	\begin{subfigure}[T]{0.18\textwidth}
		\centering
		\includegraphics[width=\textwidth]{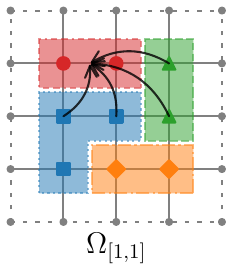}
		\caption{Spatial Communication}
	\end{subfigure}%
	\hfill
	\begin{subfigure}[T]{0.36\textwidth}
		\centering
		\includegraphics[width=\textwidth]{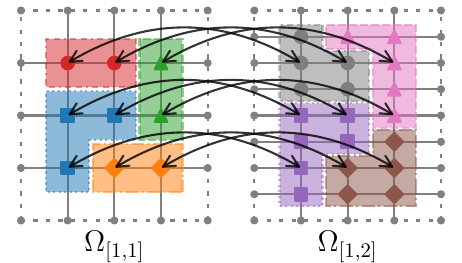}
		\caption{Subproblem Communication}
	\end{subfigure}%
	\hfill
	\begin{subfigure}[T]{0.36\textwidth}
		\centering
		\includegraphics[width=\textwidth]{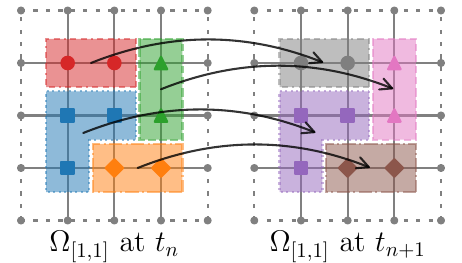}
		\caption{Temporal Communication}
	\end{subfigure}
    \caption{The communication pattern for each parallelization scale. The process association of each grid point is given by color.}\label{fig:communication_pattern}
\end{figure}

Overall, by combining the parallel-in-time method, the sparse grid combination method, and the domain decomposition method, we obtain a parallel solution method for parabolic problems \cref{pdgl} that allows three scales of nested parallelization, since both  spatial and temporal dimensions are treated in parallel. This opens a way to efficiently use an extremely large number of cores. Moreover, the associated communication pattern for the different parallelization scales is characterized as follows: For spatial parallelization, communication is local to each subproblem, since only processes involved in the same subproblem need to exchange data. For 
the combination method, communication is global, but only on the subproblem scale, i.e.\ all subproblems have to communicate with each other, but not all processes of one subproblem have to communicate with all processes of another subproblem. Communication on the time scale is again local to each subproblem. The result is an extremely parallel solver with excellent speed-up and scale-up properties that is perfectly suited for modern large-scale massively parallel computing systems of the exascale age. It allows solving Fokker-Planck problems for several important practical applications, bypassing parallel stochastic methods with much inferior accuracy and parallelization properties. The communication patterns for the different parallelization scales are depicted in \cref{fig:communication_pattern} for exemplary subproblems of the combination method.

The remainder of this paper is organized as follows: 
In \cref{sec:framework} we give the theoretical framework and the details of the building blocks of our new approach.
In \cref{sec:pitct} we present the parallel-in-time combination method for two algorithmic variants, the multigrid reduction-in-time on the sparse grid and the multigrid reduction-in-time on the subproblems.
In \cref{sec:numerical_experiments} we give the results of our parallelization experiments. We consider the heat equation, the chemical master equation, and certain exemplary stochastic differential equations.
In \cref{sec:concluding_remarks} we give some concluding remarks.

\section{Theoretical framework}
\label{sec:framework}
We consider a general parabolic problem of the form
\begin{equation}\label{eq:parabolic_problem}
\begin{gathered}
	\frac{\partial u}{\partial t} + {\cal L}u = f \, \text{in } \Omega\times(T_\text{start}, T_\text{end}]\\
    u = \bar{u}_{\Gamma_D} \, \text{on } \partial\Omega\times(T_\text{start}, T_\text{end}],
    u = \bar{u}_0 \, \text{in } \Omega\times\{T_\text{start}\}
\end{gathered}
\end{equation}
on some domain $\Omega \subset \mathbb{R}^d$ for start and end time $0 \leq T_\text{start} < T_\text{end}$, some linear elliptic differential operator $\cal L$, some Dirichlet boundary conditions $\bar{u}_{\Gamma_D}$ and initial condition $\bar{u}_0$. We are particularily interested in the case $d>3$ focusing on moderate values of $d$ up to $d=6$ which precludes any conventional space discretization on uniform grids due to the curse of dimension.

We discretize \cref{eq:parabolic_problem} for some grid $\Omega_l$ on $\Omega$, with mesh size $h_l\approx (2^{-l_1},\ldots, 2^{-l_d})$, where $l=(l_1, \ldots, l_d)$ is a $d$-dimensional multi-index parameterizing the spatial discretization. This results in a semi-discrete problem, namely the coupled system of ordinary differential equations
\begin{equation}\label{eq:parabolic_problem_discretized}
    \begin{gathered}
	    \frac{\partial u_l}{\partial t} + {\cal L}_lu_l = f_l \, \text{in } \Omega_l\times(T_\text{start}, T_\text{end}]\\
        u_l = \bar{u}_{\Gamma_D} \, \text{on } \partial\Omega_l\times(T_\text{start}, T_\text{end}],
        u_l = \bar{u}_0 \, \text{in } \Omega_l\times\{T_\text{start}\}.
    \end{gathered}
\end{equation}
Next, let $T_\text{start} = t_{l, 0} < t_{l, 1} < \ldots < t_{l, N_l} = T_\text{end}$ be a partition of the time interval of interest into $N_l + 1$ time steps $t_{l, n}, 0 \leq n \leq N_l$ with time step size $\Delta t_{l, n} \coloneq T_{l, n} - T_{l, n-1}, 1 \leq n \leq N_l$. Note that the time partition may depend on the spatial discretization parameter $l$.
Now denote by $\Phi_{l, n}$ the time propagator for time step $n$ of the semi-discrete problem \cref{eq:parabolic_problem_discretized}, i.e.\
\begin{equation}\label{eq:parabolic_problem_time_propagation}
    \begin{aligned}
    u_{l, 0} &= g_{l, 0} \coloneq \bar{u}_0 \\
    u_{l, n} &= \Phi_{l, n}(u_{l, n - 1}) + g_{l, n}, \quad 1 \leq n \leq N
    \end{aligned}
\end{equation}
for some forcing term $g_{l, n}$. In this paper we will mainly use the backward Euler method, where we have $\Phi_{l, n} = (1 + \Delta t_{l, n} {\cal L}_{l, n})^{-1}$ and $g_{l, n} = (1 + \Delta t_{l, n} {\cal L}_{l, n})^{-1} \Delta t_{l, n} f_{l, n}$. However, other linear time propagators could be used analogously. Thus, at each time step $n$ we encounter the linear system of equations
\begin{equation}\label{eq:parabolic_problem_time_propagation_LES}
    \begin{aligned}
\tilde {\cal L}_{l, n} u_{l, n} := 
	    ({\cal L}_{l, n}+\frac 1 {\Delta t_{l, n}} {\cal I}_{l,n} ) u_{l, n} &= f_{l, n} +\frac 1 {\Delta t_{l, n}}  u_{l, n-1} =: \tilde f_{l, n} \quad 1 \leq n \leq N 
    \end{aligned}
\end{equation}
with the matrices ${\cal L}_{l, n}$, ${\cal I}_{l, n}$ and the vectors $f_{l, n}$ of size about $2^{l_1+\ldots+l_d}$ discretizing the operator $\cal L$, the identity $\cal I$, and the right-hand-side function $f$ in space on grid $\Omega_l$ at time point $n$. 

\subsection{Combination method in space for elliptic problems and domain decomposition Solver}
\label{sec:ct_elliptic}
\subsubsection{Fundamentals of the combination method} 
\label{sec:ct_fundamentals}
The following is a brief summary of the combination method. 
For details, the reader is referred to \cite{Griebel2023} and the references cited therein.
The combination method provides a solution to elliptic sparse grid problems in $d$ dimensions by combining solutions to subproblems on (in general) anisotropic grids.
In particular, consider a domain $\Omega\subset\mathbb{R}^d$ and let $l=(l_1,\ldots,l_d)\in\mathbb{N}_+^d$ be a level parameterizing a collection of grids $\Omega_l\subset\Omega$ with grid sizes $h=(h_1,\ldots,h_d), h_j\approx2^{-l_j},j=1,\ldots,d$. Then the subproblem associated with $l$ consists of solving the elliptic problem of interest \cref{eq:parabolic_problem_discretized} for the solution $u_l$ on the subproblem grid $\Omega_l$.
With the usual $\ell_1$ norm $\norm{l}_1 \coloneq l_1+\ldots+l_d$ we get the approximate sparse grid solution $u_L^{(c)}$ by combining all subproblem solutions $u_l$ by the combination method according to the formula
\begin{equation}\label{eq:combination_technique_solution}
u(x) \approx u_L^{(c)}(x) \coloneqq \sum_{w=0}^{d-1} (-1)^w \binom{d-1}{w} \sum_{\substack{\norm{l}_1=L+(d-1)-w\\l \geq L_0}} u_l(x)
\end{equation}
for some level $L\in\mathbb{N}$ and initial level $L_0 \in \mathbb{N}$, $1 \leq L_0 \leq L$, where $l \geq L_0$ is to be understood component-wise. If not stated otherwise, we assume that $L_0 = 1$, i.e.\ we consider the complete selection of subproblems.
In some cases it is necessary to discard extremely anisotropic grids by choosing $L_0 > 1$, e.g.\ to be able to adequately represent initial conditions on each subproblem grid $\Omega_l$. \cref{fig:combination_technique} shows an example construction of $u_L^{(c)}$. 

Note that all subproblems are independent of each other and can be computed completely in parallel. The cost of the combination method \cref{eq:combination_technique_solution} must be compared with that of a conventional discretization of a grid with uniform mesh size $2^{-L}$ for each coordinate direction. There, the number of degrees of freedom scales as $O(2^{dL})$, which expresses the well-known curse of dimensionality that makes such a conventional discretization practically impossible for $d>3$.
In contrast, each of the grids involved in the combination method has only $O(2^L)$ degrees of freedom. Moreover, the number of involved different grids $\Omega_l$ is $O(L^{d-1})$, and each of the associated problems can be computed independently in parallel. Only the combination \cref{eq:combination_technique_solution} of the computed different solutions $u_l$ requires some sequential computations and some communication. The combination method has been shown in \cite{Griebel.Harbrecht:2014} to provide an approximation for the solution of elliptic problems that has the same convergence rate and the same order of error as the corresponding sparse grid approximation by the Galerkin method. Furthermore, it is known that the Galerkin sparse grid approximation with piecewise linear basis functions leads to an error of $O(2^{-2L} L^{d-1})$ with respect to the $L_2$ norm, provided that the solution of the considered continuous elliptic PDE is in the Sobolev space of bounded 2nd mixed deriatives ${\cal H}^2_{mix}$, see \cite{Griebel2023}. For our application, this is indeed the case (for smooth boundary conditions, and perhaps except for the very first time step $n=1$ for non-smooth initial conditions). This is due to the smoothing property of the elliptic operator $\cal L$ over time for our parabolic problem \cref{eq:parabolic_problem}. This convergence rate must be compared to the rate $O(2^{-2L})$, which would be obtained for a conventional discretization with piecewise linear basis functions on a grid of uniform mesh size $(2^{-L},\ldots,2^{-L})$. We see that the combination method loses only marginally in terms of error rate (i.e.\, by a factor of $L^{d-1}$), but gains enormously in terms of degrees of freedom involved, i.e.\ it avoids the curse of dimension in the leading cost term. In this paper, we restrict ourselves to the case of piecewise linear basis functions for simplicity. Note however that there are sparse grid combination methods with convergence rates of higher order $s$ based on piecewise polynomials of higher degree for smoother solutions in spaces ${\cal H}^s_{mix}$, $s>2$.

\subsubsection{Spatial subproblem solver by domain decomposition}
\label{sec:spatial_subproblem_solver}

For our application \cref{eq:parabolic_problem_time_propagation_LES}, all subproblems with multi-index $l$ in the above combination method involve, for any fixed time step $n$, the solution of an elliptic problem discretized on the grid $\Omega_l$, i.e.\ each subproblem~$l$ can be understood as solving 
\begin{equation}
\label{eq:ct_subproblem_discretized}
	\tilde{\cal L}_{l,n}u_{l,n} = \tilde{f}_{l,n}
\end{equation}
for $u_{l,n}$ with the given discretized operator $\tilde{\cal L}_{l,n}$ and the discretized right-hand side $\tilde{f}_{l,n}$ for time step $n$.
Thus, we need a robust and efficient solver for the large number of subproblems on generally anisotropic grids involved in the combination method \cref{eq:combination_technique_solution}, i.e.\ for the corresponding linear systems of equations \cref{eq:ct_subproblem_discretized}.
In this work we use our dimension-oblivious domain decomposition method presented in \cite{Griebel2023} as the preconditioner for a Krylov iterative solution method. The solver, equipped with our preconditioner, is capable of efficiently handling the solution procedure on all subproblem grids generated by the combination method. Alternatively, other methods such as multigrid or BPX-multilevel can be used. However, such schemes are more complicated to parallelize efficiently for arbitrary $\Omega_l$ and in particular on large compute systems compared to our domain decomposition approach. Additionally, their achieved convergence rate is, up to the order constant, not better than for our domain decomposition scheme.
In the following we give a short overview of the preconditioning algorithm, for details see \cite{Griebel2023}.

Let $\Omega_l$ be any, usually anisotropic, subproblem grid in arbitrary dimension $d$ that appears in the combination method, and let $\tilde{L}_{l,n}$ and $\tilde{f}_{l,n}$ be the corresponding discretized elliptic operator and right-hand side at time point $n$.
We first generate an initial partition $\{\bar{\Omega}_{l, i}\}_{i = 1}^{P^x_l}$ of $\Omega_l$ by dividing the grid points along a discrete space-filling curve, i.e., the Hilbert curve in $d$ dimensions, into $P^x_l \geq 1$ (approximately) equal-sized disjoint subdomains $\bar{\Omega}_{l, i}$ of size $\varpi_l := \lfloor\frac{\lvert\Omega_l\rvert}{P^x_l}\rfloor$
 or $\varpi_l + 1$.
In fact, we obtain exactly 
$P_l^x - (\lvert{\Omega_l}\rvert - P_l^x \varpi_l)$ 
subdomains of size 
$\varpi_l$
and 
$\lvert\Omega_l\rvert- P^x_l \varpi_l]$ 
subdomains of size 
$\varpi_l + 1$.

\begin{figure}[htb]
    \centering
    \includegraphics[]{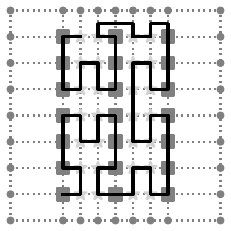}\hfill
    \includegraphics[]{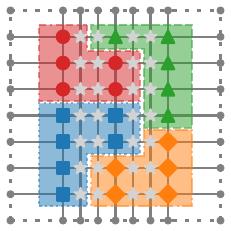}\hfill
    \includegraphics[]{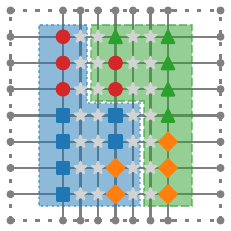}\hfill
  \caption{Decomposition of an anisotropic grid with $l=(2,3)$ by the Hilbert curve approach (left: ordering of nodes (including inner virtual nodes) by the space-filling curve, center: disjoint partitioning of the nodes (color-coded by owner process) and the respective processor subdomains (color-coded by owner process) by splitting the curve into four intervals). Enlargement of these disjoint subdomains by $\gamma = 0.25$ and the resulting overlapping subdomains (right: only two of the four subdomains are color-coded for visibility, all nodes however are still color-coded by owner). Virtual inner nodes are depicted by stars.}
    \label{fig:ddm_partition}
\end{figure}

Then, each of these subdomains is enlarged by an overlap factor $\gamma > 0$, which controls how many grid points of the neighboring disjoint subdomains $\bar{\Omega}_{l, i \pm k}, k=1, \ldots, \lceil 2\gamma \rceil$ along the space-filling curve must be added to the current subdomain $\bar{\Omega}_{l, i}$ to generate the extended subdomain $\Omega_{l, i}$. The final overlapping domain decomposition of $\Omega_l$ is then given by the collection of these $\Omega_{l, i}$.
This construction is illustrated for $d=2$ in \cref{fig:ddm_partition}. 

Note that the Hilbert curve is isotropic, but we encounter anisotropic grids in the combination method. To deal with the anisotropy of a grid $\Omega_l$, we adjust in this paper the non-maximum dimensions to the maximum dimension by inserting an appropriate number of {\em virtual} nodes uniformly between the actual grid nodes, which results in an isotropic grid. Note that these virtual nodes are only used for the enumeration of the actual grid nodes along the space-filling curve and are not inserted into the grid in praxis. Thus they can be interpreted as an offset in the enumeration along the space-filling curve. This is in contrast to \cite{Griebel2023}, where the anisotropic grids were padded by padding each non-maximum dimensions from the right to match the maximum dimension. Our new approach allows for more uniformly shaped subdomains, which slightly improves the convergence rate of the spatial solver.
The specific choice
\begin{equation}
    \label{eq:gamma_choice}
    \gamma = \frac{1}{2}m
\end{equation}
for some $m\in\mathbb{N}, m\leq P^x_l - 1$, ensures that each grid point is contained in exactly $m+1$ extended subdomains, see \cite{Griebel2023}. This also allows the solver to recover from hardware failure of $m$ neighboring processes in the same iteration, which is an ongoing research topic. The restriction $R_{l, i} \in \mathbb{R}^{\lvert \Omega_{l, i} \rvert \times \lvert \Omega_l \rvert}$ and prolongation $P_{l, i} \coloneq R_{l, i}^T \in \mathbb{R}^{\lvert \Omega_l \rvert \times \lvert \Omega_{l, i} \rvert}$ corresponding to subdomain $\Omega_{l, i}$, $i = 1, \ldots, P_l^x$, are simply given by the trivial injection operators.

Besides a decomposition of the domain of the spatial subproblem on the fine scale like the partition $\{ \Omega_{l, i}, i = 1, \ldots,{P^x_l}\}$, a good domain decomposition method also needs a {\em coarse space problem} for fast convergence. 
With $q_l$ denoting the number of coarse degrees of freedom per fine-scale subdomain so that the coarse space contains $q_l\cdot P^x_l$ degrees of freedom, the algebraic restriction to the coarse space for subproblem $l$ is given by
\begin{equation}
R_{l, 0} \coloneq \blockdiag_{i=1}^{P_l^x} (R_{l, i, 0}) \in \mathbb{R}^{(P_l^x \cdot q_l) \times (\prod_{i=1}^{P_l^x} \lvert \Omega_{l, i} \rvert)},
\end{equation}
where $R_{l, i, 0} \in \mathbb{R}^{q_l \times \lvert \Omega_{l, i} \rvert}$ and
\begin{equation}
R_{l, i, 0} \coloneq
\begin{pNiceMatrix}[first-row]%
    \rule[-8pt]{0pt}{8pt}
    \Block{1-10}{q_l\text{ blocks of size }\approx \lvert \bar{\Omega}_{l, i} \rvert/q_l} \\
    1 &  \ldots & 1 & \phantom{1} & & \phantom{1} & & & & \phantom{1} \\
      &         & \phantom{1}  & 1 & \ldots & 1 \\
      &         &   &  &  & &  \ddots \\
      &         &   &  &  &  & & 1 & \ldots & 1
    \CodeAfter
       \begin{tikzpicture}
	       \draw[brace, line width=0.5pt]
            ([xshift=3mm]1-10.north east) to node[auto = left] {$~q_l$ rows}
            ([xshift=3mm,yshift=-1.4mm]4-10.south east);
       \draw[brace=4pt, line width=0.5pt]
            (1-1.north west) to node[above=3pt] {}
            (1-10.north east);
       \end{tikzpicture}
    \end{pNiceMatrix}\qquad\qquad.
\end{equation}
Note that we associated here the coarse problem to the index $i = 0$ and we associated the subdomain problems to the indices $i=1, \ldots, P_l^x$. Next, we define the local operators as $\tilde{\cal L}_{l, n, i} \coloneq R_{l, i} \tilde{\cal L}_{l,n} R_{l, i}^T, i = 0, \ldots, P_l^x$.
Let $C_{l,n, (1), D}^{-1}$ be the one-level preconditioner given by
\begin{equation}
C_{l,n, (1), D}^{-1} \coloneq \sum_{i = 1}^{P_l^x} R_{l, i}^T D_{l, i} \tilde{\cal L}_{l,n, i}^{-1}R_{l, i},
\end{equation}
where $D_{l, i}$ are diagonal matrices $D_{l, i}\in\mathbb{R}^{N_{l, i}\times N_{l, i}}$ dealing with the multiplicity of degrees of freedom due to the overlapping nature of the subdomains.  
In this paper, we restrict ourselves to $D_{l,n, i} = \omega_{l, i} I$, $i=1, \ldots, P_l^x$, where $\omega_{l, i} = 1/\max_{j\in\bar{\Omega}_{l, i}}
(\lvert\{\Omega_{l, k}, k=1,\ldots,P^x_l: j\in\bar{\Omega}_{l, k}\}\rvert)$, i.e.\ $\omega_{l, i}$ is the inverse of the maximum number of subdomains that overlap any point contained in $\bar{\Omega}_{l, i}$.
Moreover, our specific choice of $\gamma$ from \cref{eq:gamma_choice} allows us to simplify to $D_{l,n, i} = \frac{1}{n+1} I$ since $\lvert\{\Omega_{l, k}, k=1,\ldots,P^x_l: j\in\bar{\Omega}_{l, k}\}\rvert \equiv n + 1$ for all grid points $j$.
In \cite{Griebel2023} other weightings $\omega_{l,i}$ were also presented and studied. 

Now let $F_{l,n} \coloneq R_{l, 0}^T \tilde{\cal L}_{l,n, 0}^{-1} R_{l, 0}$. Then our additive two-level domain decomposition preconditioner is given by
\begin{equation}
	C_{l,n, (2), D, add}^{-1} \coloneq C_{l,n, (1), D}^{-1} + F_{l,n}.
\end{equation}
Finally, let $G_{l,n} \coloneq I - \tilde{\cal L}_{l,n}F_{l,n}$. Then a balanced version of our two-level domain decomposition preconditioner, e.g.\ with further improved convergence properties, is
\begin{equation}
    \label{eq:preconditioner_balanced}
	C_{l,n, (2), D, bal}^{-1} \coloneq G_{l,n}^TC_{l,n, (1), D}^{-1}G_{l,n} + F_{l,n}.
\end{equation}
Note that the choice of the number of subdomains $P^x_l$, and thus the number of processes used, is crucial for the parallel efficiency of the spatial solver when applied to each subproblem $l$ of the combination method.
Recall that all subproblems are solved simultaneously. For a uniform distribution of the total load, each process should get an equal amount of work. However, the number of degrees of freedom varies between different subproblem layers, i.e.\ different values $w$ in the combination method \cref{eq:combination_technique_solution} (but not within each layer). In \cite{Griebel2023} it was shown that
\begin{equation}
    \label{eq:spatial_speedup_factor}
    P^x_l \coloneq \hat{P}^x \cdot 2^{d - w - 1}
\end{equation}
for some spatial parallelization factor $\hat{P}^x\in\mathbb{N}_{>1}$ results in disjoint subdomain sizes $\lvert\bar{\Omega}_{l, i}\rvert \approx \frac{2^L}{\hat{P}^x}$, so the enlarged subdomain will contain $\lvert\Omega_{l, i}\rvert \approx (1 + 2\gamma)\frac{2^L}{\hat{P}^x}$ nodes.
Note that with this particular choice, the size of the enlarged subdomains of each subproblem, i.e.\ the number of nodes stored on each process, is independent of dimension $d$ and layer $w$, and in particular independent of subproblem~$l$, leading to an approximately balanced use of all processes.
The results of \cite{Griebel2023} show that this gives a robust, dimension-independent preconditioner with optimal order convergence behavior and excellent scalability properties.
Finally note that the number of processes $P^x_l$ of \cref{eq:spatial_speedup_factor} is constant for each subproblem layer $w$, i.e.\ it does not depend on the subproblem $l$. However, since the grids of the (extremely) anisotropic subproblems on each layer $w$ contain fewer grid points than the more isotropic grids on the same layer due to boundary handling,
the number of grid points assigned to each spatial process can still differ. For example, consider the subproblems on layer $w = 0$ for $d = 2$ and $L = 11$. One of the extremely anisotropic grids is associated with the subproblem $l = (1, 11)$. It contains one interior node along the first dimension and $2.047$ interior nodes along the second dimension, for a total of $2.047$ nodes.
In contrast, the isotropic grid $l = (6, 6)$ contains $63$ interior grid nodes in each dimension, for a total of $3.969$ nodes.
However, for the domain decomposition of problems $\hat{P}^x \cdot 2^{d-w-1}=\hat{P}^x \cdot 2$
spatial processes are used in parallel, since they belong to the same layer. Therefore, the subdomains of the isotropic subproblem will be almost twice as large as those of the anisotropic subproblem.
To improve this aspect, we use a subproblem-dependent choice of the number of processes. For this purpose, denote by $2^S \leq \lvert \Omega_l \rvert$ the target number of nodes per subdomain, i.e.\ we want the disjoint subdomains to be of size $\lvert\bar{\Omega}_{l, i}\rvert \approx 2^S$.
To achieve this, we choose
\begin{equation}\label{eq:spatial_speedup_factor_new}
	P^x_l \coloneq \lceil\frac{\lvert\Omega_l\rvert}{2^S}\rceil = \lceil\frac{\prod_{j=1}^{d}(2^{l_j} - 1)}{2^S}\rceil,
\end{equation}
so that $\lvert\bar{\Omega}_{l, i}\rvert = \frac{\lvert\Omega_l\rvert}{P^x_l} = \frac{\lvert\Omega_l\rvert}{\lceil\frac{\lvert\Omega_l\rvert}{2^S}\rceil} \approx 2^S$.
Looking again at the two subproblems above and choosing the target number of nodes as $2^{10}$, i.e.\ $S = 10$, the subproblem $l = (1, 11)$ will now use two spatial processes, while $l = (6, 6)$ will now use four spatial processes, resulting in a balanced subdomain size of about $1.024$ nodes for each subproblem.

\subsection{Time-sequential combination method}
\label{sec:sequential_combination_technique}
In the following we describe a general approach for the time-sequential solution of parabolic problems using the combination method given in \cref{sec:ct_fundamentals} and the spatial subproblem solver discussed in \cref{sec:ct_elliptic} for the discretized elliptic problems \cref{eq:ct_subproblem_discretized}.

\subsubsection{General algorithm}
\label{sec:general_algorithm}
\begin{figure}[htb]
    \centering
    \includegraphics[]{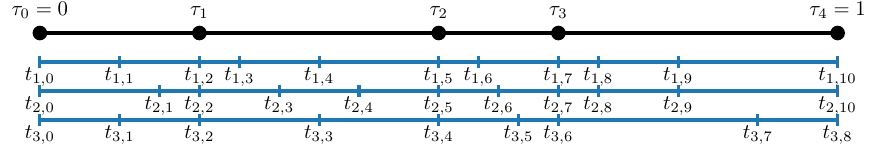}
    \caption{Partition of the time interval $[0, 1]$ into recombination steps and an example of subproblem-dependent time steps for three subproblems.}
    \label{fig:ct_time_partition}
\end{figure}

Let $T_\text{start} = \tau_0 < \tau_1 < \ldots < \tau_s = T_\text{end}, s \geq 1$, be a partition of the time interval of interest $[T_\text{start}, T_\text{end}]$ into $s + 1$ so-called recombination steps $\tau_k$.
Following the notation from the beginning of this section, let $T_\text{start} = t_{l, 0} < t_{l, 1} < \ldots < T_{l, N_l} = T_\text{end}$ be a partition of $[T_\text{start}, T_\text{end}]$ independently for each subproblem $l$ of the combination method.
Denote by $0 \leq n(l, k) \leq N_l$ the unique time step index such that $t_{l, n(l, k)} = \tau_k, 0 \leq k \leq s$. We require that $n(l, k)$ is well defined for all $k$ and all subproblems $l$. This ensures that all recombination steps $\tau_k$ are present in the time partitions of all subproblems. In particular we have $n(l, 0) = 0$ and $n(l, s) = N_l$. An example of time partitions for $M=3$ subproblems, here for simplicity indexed as $l=1,2,3$, is illustrated in \cref{fig:ct_time_partition}.

\begin{figure}[htb]
    \centering
    \includegraphics[]{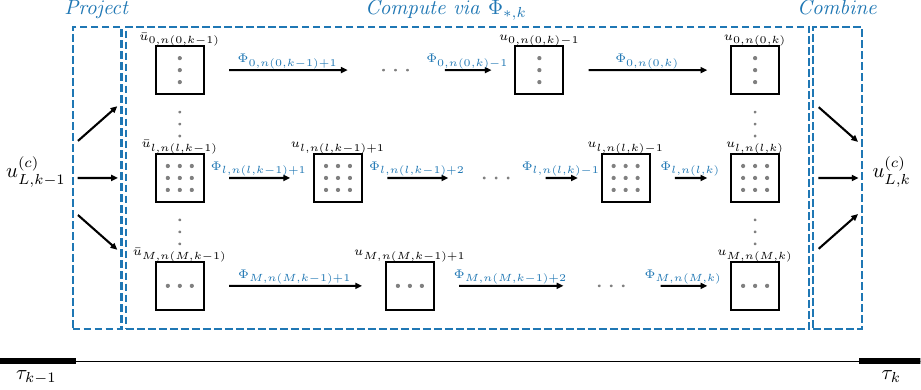}
    \caption{Overview of one time propagation step of the time-sequential combination method in \cref{alg:combination_technique_general_parabolic} for various subproblems.}
    \label{fig:combination_technique_general_parabolic}
\end{figure}

In the exemplary time partition in \cref{fig:ct_time_partition} the time step for the subproblem $l=1$ corresponding to the recombination step $\tau_2$ is $t_{1, n(1, 2)} = t_{1, 5}$. Note that by construction the time steps corresponding to recombination steps have to coincide and time steps between recombination steps can be chosen arbitrarily and independently for each subproblem, allowing for different time stepping schemes for each subproblem.
Let $\Phi_{l, n}$ be the time propagators for the considered parabolic problem with time steps $t_{l, n}, 1 \leq n \leq N_l,$ for subproblem $l$ analogous to \cref{eq:parabolic_problem_time_propagation}.
Using this notation, the solution of \cref{eq:parabolic_problem} by the combination method in space in a time-sequential manner is given by \cref{alg:combination_technique_general_parabolic}. A graphical illustration is given in \cref{fig:combination_technique_general_parabolic}.
\begin{algorithm}[hbt]
    \caption{Time-Sequential Combination method for Parabolic Problems}
    \label{alg:combination_technique_general_parabolic}
    \hspace*{\algorithmicindent} \textbf{Input} Spatial subproblem grids $\Omega_l$ for the combination method with sparse grid level $L$\\
    \hspace*{\algorithmicindent} \hphantom{\textbf{Input}} Time propagators $\Phi_{l, n}$ and time steps $t_{l, n}, 1 \leq n \leq N_l$, for each subproblem $l$\\
    \hspace*{\algorithmicindent} \hphantom{\textbf{Input}} Recombination steps $\tau_k$, $0 \leq k \leq s$ \\
    \hspace*{\algorithmicindent} \hphantom{\textbf{Input}} Initial value $\bar{u}_0$ \\
    \hspace*{\algorithmicindent} \textbf{Output} Approximate sparse grid solution $u^{(c)}_{L, s}$ at recombination step $\tau_s \coloneq T_\text{end}$
    \begin{algorithmic}[1]
        \State Set $u_{l, n(l, 0)}$ to the initial value $\bar{u}_0$ for all subproblems $l$.
	    \State Combine all $u_{l, n(l, 0)}$ using \cref{eq:combination_technique_solution} to obtain $u^{(c)}_{L, 0}$. \Comment{inter-subproblem comm.}\label{alg:ct_general_parabolic_combine_initial}
        \For{$1 \leq k \leq s$}
            \For{each subproblem $l$ simultaneously}
                \State Project $u^{(c)}_{L, k - 1}$ onto $\Omega_l$ to obtain $\bar{u}_{l, n(l, k - 1)}$. \label{alg:ct_general_parabolic_project}
                \State Compute $u_{l, n(l, k)}$ from $\bar{u}_{l, n(l, k-1)}$ via the time propagators $\{\Phi_{l, n(l, k - 1) + 1}, \ldots, \Phi_{l, n(l, k)}\}$.
            \EndFor
	    \State Combine all $u_{l, n(l, k)}$ using \cref{eq:combination_technique_solution} to obtain $u^{(c)}_{L, k}$. \Comment{inter-subproblem comm.}\label{alg:ct_general_parabolic_combine}
        \EndFor
        \State\Return $u^{(c)}_{L, s}$
    \end{algorithmic}
\end{algorithm}

\subsubsection{Combining spatial subproblem solutions}
\label{sec:hierarchisation}
The missing ingredient for the implementation of \cref{alg:combination_technique_general_parabolic} is an efficient combination of all intermediate subproblem solutions $u_{l, n(l, k)}$ to obtain $u^{(c)}_{L, k}$, as required in \cref{alg:ct_general_parabolic_combine_initial} and \cref{alg:ct_general_parabolic_combine} of \cref{alg:combination_technique_general_parabolic}, and projections of $u^{(c)}_{L, k - 1}$ onto the subproblem grids $\Omega_l$ resulting in $\bar{u}_{l, n(l, k - 1)}$, as required in \cref{alg:ct_general_parabolic_project} of \cref{alg:combination_technique_general_parabolic}.
Note that the combination of the subproblem solutions is immediately followed by a projection step from the combination solution to each subproblem in all iterations. By using a basis transformation in space to a hierarchical basis representation of the coefficients on the respective grids, we can perform the combination and the projection at once, which allows to avoid the construction of the intermediate solutions  $u^{(c)}_{L, k}$ for $1 \leq k < s$ of the combination method. In the following we give a short overview, details can be found in \cite{Garcke:2013,Griebel1992}.

To simplify the notation we drop the time indices and assume homogeneous Dirichlet boundary conditions,\footnote{Inhomogeneous boundary conditions can be treated in a similar way.} so that we can consider only interior basis functions. Let $V_l$ denote the space of piecewise d-linear interior nodal basis functions associated with the subproblem grid $\Omega_l$. Then we define by
\begin{equation}
    W_l \coloneq V_l \setminus \oplus_{j=1}^d V_{l - e_j}
\end{equation}
the hierarchical difference space $W_l$, where $e_j$ is the unit vector in the direction of the $j$-th dimension. We set $V_l \coloneq 0$ if any component of $l$ is negative.
\begin{figure}
	\centering
    \includegraphics[]{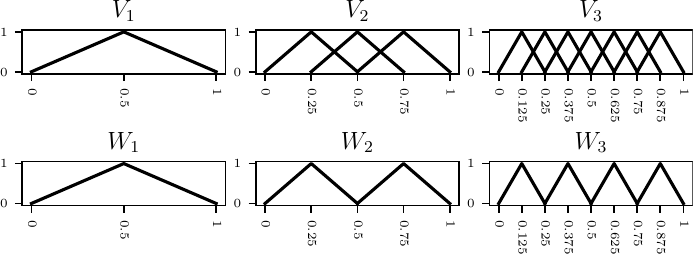}
    \caption{One-dimensional nodal and hierarchical basis functions for levels $l=1,2,3$.}
    \label{fig:hierarchical_basis}
\end{figure}
\cref{fig:hierarchical_basis} shows a one-dimensional nodal basis for $V_l$ and associated hierarchical basis for $W_l$, $l=1,2,3$. Note that the size of the support of each nodal basis function associated with each node of $V_l$ is $2^{-l+1}$ and changes with $l$, while it is $2^{-l'+1}$ for the functions in $W_{l'}, l'=1,2,l$, which altogether form the hierarchical basis. Since $V_l$ is the union of all hierarchical difference spaces $W_{l'}$ such that $l' \leq l$ componentwise it can be associated with a hierarchical basis. Due to this we have that the hierarchical basis functions associated to coarser levels of $V$ are simply a subset of the hierarchical basis functions associated to finer levels of $V$.
This allows to replace the combination of subproblem solutions on the function level by a simple combination of the coefficients of their representations in the hierarchical basis via basically the same formula \cref{eq:combination_technique_solution}, i.e.\ just the coefficient values of the representation in the hierarchical basis of a discrete function must be properly added or subtracted.

In particular, after solving the subproblem on $\Omega_l$ in the nodal basis $V_l$ and obtaining $u_l$, we hierarchize its coefficients $\tilde{u}_l\in\mathbb{R}^{\lvert\Omega_l\rvert}$ associated with the nodal basis to obtain coefficients $\hat{u}_l\in\mathbb{R}^{\lvert\Omega_l\rvert}$ associated with the hierarchical basis. These hierarchical coefficients can then be combined using \cref{eq:combination_technique_solution} to obtain the hierarchical coefficients $\hat{\bar{u}}_l\in\mathbb{R}^{\lvert\Omega_l\rvert}$ of $\bar{u}_l$, the projection of $u^{(c)}_L$ onto $\Omega_l$. A de-hierarchization of $\hat{\bar{u}}_l$ then yields the corresponding nodal coefficients $\tilde{\bar{u}}_l\in\mathbb{R}^{\lvert\Omega_l\rvert}$. The processes of hierarchization and de-hierarchization thus involves basically the linear transformations between the representations in the conventional nodal basis and the hierarchical basis and are given in detail in e.g.\ \cite{Griebel:1998}.
Note that they require communication between all subproblems, since, for example, the central grid point is contained in each subproblem grid. However, in view of the parallel spatial subproblem solver of \cref{sec:spatial_subproblem_solver}, only those processes of each subproblem that contain the same node need to exchange data, i.e.\ each process of each subproblem only needs to communicate with a subset of processes of other subproblems. More details can be found in \cite{Rentrop.Griebel.2020}.

Now denote by $\tilde{u}^{(c)}_L \coloneq \blockdiag_l (\tilde{u}_l) \in \mathbb{R}^{\prod_l\lvert\Omega_l\rvert}$ the collection of all nodal coefficients $\tilde{u}_l$ of the subproblems, analogously for $\hat{u}^{(c)}_L$. Let $H^{(c)}_l \in \mathbb{R}^{\lvert \Omega_l \rvert \times \lvert \Omega_l \rvert}$ and $H^{(c), -1}_l \in \mathbb{R}^{\lvert \Omega_l \rvert \times \lvert \Omega_l \rvert}$ be the local hierarchization and de-hierarchization operators on $\Omega_l$, i.e.\ $\hat{u}_l = H^{(c)}_l\tilde{u}_l$ and $\tilde{u}_l = H^{(c), -1}_l\hat{u}_l$. Furthermore, let $H^{(c)}_L = \blockdiag_l (H^{(c)}_l) \in \mathbb{R}^{(\prod_l\lvert\Omega_l\rvert) \times (\prod_l\lvert\Omega_l\rvert)} $ and $H^{(c), -1}_L = \blockdiag_l (H^{(c), -1}_l) \in \mathbb{R}^{(\prod_l\lvert\Omega_l\rvert) \times (\prod_l\lvert\Omega_l\rvert)}$ denote the global hierarchization and global de-hierarchization operators such that $\hat{u}^{(c)}_L = H^{(c)}_L \tilde{u}^{(c)}_L$ and $\tilde{u}^{(c)}_L = H^{(c), -1}_L \hat{u}^{(c)}_L$.
Next, denote by $S(i, \Omega_l) \coloneq \{(j, \Omega_m): j \in \Omega_m, x_j = x_i\}$ the nodes on all subproblem grids that correspond to the node $i$ on $\Omega_l$, where $x_i \in \mathbb{R}^d$ denotes the coordinates of the $i$-th node in $\Omega_l$, and let $\mathrm{layer}(\Omega_l) \coloneq L + (d - 1) - \norm{l}_1$ be the layer of subproblem $l$.
Finally, let $\bar{Q}_L$ be the combination operator representing \cref{eq:combination_technique_solution}, i.e\
\begin{equation}
    \label{eq:recombination_operator}
    (\bar{Q}_L\tilde{u}^{(c)}_L)_{(l, i)} \coloneq \sum_{(j, \Omega_m)\in S(i, \Omega_l)} (-1)^{\mathrm{layer}(\Omega_m)} \binom{d - 1}{\mathrm{layer}(\Omega_m)}(\tilde{u}^{(c)}_L)_{(m, j)}.
\end{equation}
Then, the process of combination followed by projection corresponding to \cref{alg:ct_general_parabolic_combine,alg:ct_general_parabolic_project} of \cref{alg:combination_technique_general_parabolic} is given by
\begin{equation}
    \label{eq:overall_combination_operator}
    Q^{(c)}_L \coloneq (H^{(c)})^{-1}_L \bar{Q}^{(c)}_L H^{(c)}_L.
\end{equation}
Its application to a global nodal coefficient $\tilde{\bar{u}}^{(c)}_L$ of some function $\bar{u}^{(c)}_L$, i.e.\ $\tilde{u}^{(c)}_L \coloneq Q^{(c)}_L\tilde{\bar{u}}^{(c)}_L$, can be understood as follows: We first hierarchize all subproblem coefficients $\tilde{\bar{u}}^{(c)}_l$ independently of each other, i.e.\ we form $H^{(c)}_L\tilde{\bar{u}}^{(c)}_L$, resulting in the hierarchical coefficients $\hat{\bar{u}}_l$. We then combine the hierarchical coefficients, $Q^{(c)}_L\hat{\bar{u}}_l$, resulting in the recombined hierarchical coefficients $\hat{\bar{u}}^{(c)}_l$. Finally, we de-hierarchize the subproblem coefficients independently of each other, i.e.\ we form $(H^{(c)})^{-1}_L \hat{\bar{u}}^{(c)}_l$, obtaining the nodal coefficients $\tilde{u}^{(c)}_L$ of a function $u^{(c)}_L$, where $u^{(c)}_L$ is the recombination and projection of $\bar{u}^{(c)}_L$ to each subproblem.
Now, reintroducing the time index, let $\Phi^{(c)}_{L, k} \coloneq \blockdiag_l(\Phi_{l, k})$ and $\tilde{u}_k \coloneq \blockdiag_l(\tilde{u}_{l, n(l, k)})$. With this, \cref{alg:combination_technique_general_parabolic} can equivalently be written as
\begin{equation}
    \label{eq:combination_technique_general_parabolic}
    \tilde{u}^{(c)}_{L, s} = Q^{(c)}_L\Phi^{(c)}_{L, s}Q^{(c)}_L\Phi^{(c)}_{L, s-1}Q^{(c)}_L\cdots\Phi^{(c)}_{L, 1}\tilde{u}^{(c)}_{L, 0}.
\end{equation}

Note that there are other techniques with improved properties for combining subproblem solutions, for example via the use of biorthogonal hierarchical basis functions as in \cite{Pollinger2023}, or via the use of prewavelets \cite{Griebel.Oswald:1995}.
For the present work, however, we use the standard hierarchical basis for simplicity.

\subsection{Multigrid reduction-in-time}
Parallel-in-time techniques offer significant opportunities for parallelization by solving the PDE in space-time, using independent discretizations of space and time instead of treating each time step in a time-sequential manner.
Prominent examples are the parallel-in-time preconditioners of \cite{Neumueller2019}, \texttt{Parareal} (\cite{Lions2001}), and \texttt{MGRIT} (\cite{Falgout2014}). In this work we use the multigrid reduction-in-time algorithm \texttt{MGRIT} of \cite{Falgout2014} via its implementation in the library \texttt{XBraid}\footnote{In this work we use the version of \texttt{XBraid} at git commit \texttt{4bbd644}.} (see \cite{xbraid-package}), for parallelizing time. \texttt{XBraid} was chosen because of its mature state and non-intrusive \texttt{C} implementation.
This allows for seamless integration with our domain decomposition and sparse grid software framework implemented in \texttt{C++}.

\subsubsection{Algorithmic overview}
\label{sec:xbraid_overview}

\texttt{XBraid} implements a variant of \texttt{MGRIT} based on a full approximation scheme (FAS). We will give a description of the two-level algorithm for simplicity, the general-level algorithm results by recursively applying the two-level scheme. We follow the presentation and notation of \texttt{XBraid}, \cite{xbraid-package}, where details can be found. For simplicity, we assume a uniform fine time partition $T_\text{start} = t_{l, 0} < t_{l, 1} < \ldots < t_{l, N} = T_\text{end}$ and $\Delta t_{l, n} = \Delta t_l$ for all $n$. Note that \texttt{XBraid} can handle variable and adaptive time step sizes as well. We call $t_{l, n}$ a fine time step and $\Delta t_{l, n}$ a fine time step size.

Recall the discretization of \cref{eq:parabolic_problem} on a fixed discretization $\Omega_l$ in space given by \cref{eq:parabolic_problem_discretized} and its iterative solution procedure with respect to time by linear time propagators $\Phi_{l,n}: \Omega_l \to \Omega_l$ from \cref{eq:parabolic_problem_time_propagation}, namely
\begin{equation}\label{eq:iterative_temporal_solution}
\begin{aligned}
	u_{l,0} &= g_{l,0}, \\
	u_{l,n} &= \Phi_{l,n}u_{l,n-1} + g_{l,n}, \quad 1 \leq n \leq N.
\end{aligned}
\end{equation}
This process is directly equivalent to the forward solution of the {\em space-time} system
\begin{equation}\label{eq:space_time_solution}
	B^{(N)}_l \vec{u}^{(N)}_l
\coloneq
\begin{pmatrix}
	{\cal I}_l \\
	-\Phi_{l,1} & {\cal I}_l & \\
    & & \ddots & \ddots & \\
	& & & -\Phi_{l,N} & {\cal I}_l
\end{pmatrix}
\begin{pmatrix}
	u_{l,0} \\ u_{l,1} \\ \vdots \\ u_{l,N}
\end{pmatrix}
=
\begin{pmatrix}
	g_{l,0} \\ g_{l,1} \\ \vdots \\ g_{l,N}
\end{pmatrix}
	\eqcolon \vec{g}^{(N)}_l,
\end{equation}
with the vectors $\vec{u}_l$ and $\vec{g}_l$ of length $N+1$ where each component $u_{l,n} \in \mathbb{R}^{\lvert \Omega_l \rvert}$ and $g_{l,n} \in \mathbb{R}^{\lvert \Omega_l \rvert}$ is associated to the repective grid $\Omega_l$. The system matrix $B^{(N)}_l$ is a block matrix of size $(N + 1) \times (N + 1)$, with ${\cal I}_l: \Omega_l \to \Omega_l$ denoting the identity.

Choosing a time coarsening factor $c > 1$ allows to define a uniform coarse time partition $T_\text{start} = T_{l, 0} < T_{l, 1} < \ldots < T_{l, N_\text{coarse}} = T_\text{end}$ with $N_\text{coarse} + 1$ coarse time steps and coarse time step size $\Delta T_{l, m'} = \Delta T_l = c \Delta t_l$, where $N_\text{coarse} = \frac{N}{c}$. The block operator of the space-time system for the coarse time partition of size $(N_\text{coarse} + 1) \times (N_\text{coarse} + 1)$ is then given by
\begin{equation}\label{eq:coarse_space_time_solution}
	B^{(N)}_{\text{coarse},l}
    \coloneq
    \begin{pNiceMatrix}
	    {\cal I}_l & & \\
	    -\Phi_{\text{coarse}, l, 1} & {\cal I}_l & & \\
         &  & \phantom{1} &  \\
	    &  & -\Phi_{\text{coarse},l, N_\text{coarse}} & {\cal I}_l
        \CodeAfter
        \begin{tikzpicture}
            \draw[dotted] ([xshift=3mm,yshift=-0.5mm]2-1.south) -- ([xshift=-3mm]4-3.north);
	       \draw[dotted] ([xshift=2mm,yshift=-0.5mm]2-2.south east) -- ([xshift=-2mm]4-4.north west);
       \end{tikzpicture}
    \end{pNiceMatrix},
\end{equation}
for appropriately chosen coarse linear time propagators $\Phi_{\text{coarse},l, m}, m=1, \ldots, N_\text{coarse}$.
A simple choice for $\Phi_{\text{coarse},l, m}$ are the fine-level time propagators $\Phi_{l,n}$, but using the coarse time step size $\Delta T_{l, m}$ instead of the fine partitioms $\Delta t_{l, n}$. For example, for the backward Euler scheme, the coarse operator corresponding to a fine time partition with the time step size $\Delta t_{l, n}$, i.e.\ $\Phi_{l, n} = (1 + \Delta t_{l, n} {\cal L}_{l, n})^{-1}$, would be given by $\Phi_{\text{coarse},l, m} = (1 + \Delta T_{l, m} {\cal L}_{l, m})^{-1}$

\begin{figure}[hbt]
	\centering
    \includegraphics[]{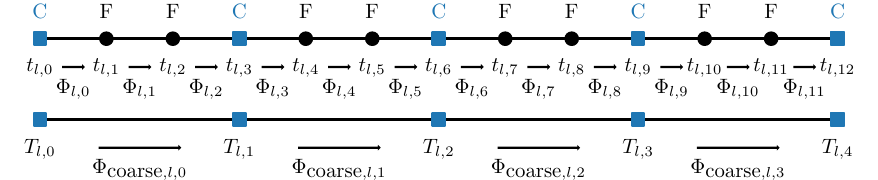}
    \caption{Classification of $N = 12$ fine time steps into F and C time steps for time coarsening factor $c = 3$.}
    \label{fig:cf_classification}
\end{figure}

Next, we split all fine time steps into $C$ time steps, i.e.\ time steps that exist in both the coarse and fine time partitions, and $F$ time steps, i.e.\ time steps that exist only in the fine time partition. This allows to define relaxation schemes and transfer operators between the two time partitions, thus completing all the necessary algorithmic ingredients for a two-level method. An example is given in \cref{fig:cf_classification}.

The components for relaxation schemes available in \texttt{XBraid} are so-called F- and C-relaxation. There, {F-relaxation} propagates the value of $u_{l,m}$ at each $C$ time step $T_m$, which is given by $u_{l, cm}$, over its corresponding coarse time interval $[T_m, T_{m+1})$ by using the sequence of fine partition time propagators ${\Phi_{l,cm + 1}, \ldots, \Phi_{l,(c + 1)m - 1}}$ successively to the values $u_{l,cm}, \ldots, u_{l,(c + 1)m - 2}$ to obtain new values of $u_l$ at all $F$ time steps $t_{cm + 1}, \ldots, t_{(c + 1)m - 1}$ contained in $[T_m, T_{m+1})$. Note that this process can be done independently for each coarse time interval, processing each interval simultaneously, but the propagation over $F$-values in each interval is inherently sequential.
C-relaxation computes the value $u_{l,cm}$ at a $C$ time step $T_m$ from the value $u_{l,cm - 1}$ at the previous $F$ time step $t_{cm - 1}$ using the fine-level time propagator $\Phi_{l,cm}$. Again, each of these updates can be done simultaneously.
Combinations of these two relaxation schemes, such as $CF$- or $FCF$-relaxation, are possible. Again, see \cite{xbraid-package} for details.

The restriction operator $R_{l, N}$ in time is given by discarding the $F$-values and the respective prolongation $P_{l, N}$ is given by injection followed by $F$-relaxation, which corresponds to harmonic interpolation in time, see \cite{xbraid-package,Falgout2014}. So we have
\begin{equation}
    R^{(N)}_l \coloneqq \begin{pmatrix}
	    {\cal I}_l \\
        0 \\
        \vdots \\
        0 \\
	    & {\cal I}_l \\
         & 0 \\
         & \vdots \\
         & 0 \\
         & & \ddots
    \end{pmatrix}^T, \mbox{ and }  
    P^{(N)}_l \coloneqq \begin{pmatrix}
	    {\cal I}_l \\
	    \Phi_{l,0} \\
	    \Phi_{l,1} \circ \Phi_{l,0} \\
        \vdots \\
	    \Phi_{l,c-2} \circ \hdots \circ \Phi_{l,0} \\
	    & {\cal I}_l \\
         & \Phi_{l,cm} \\
         & \Phi_{l,cm + 1} \circ \Phi_{l,cm} \\
         & \vdots \\
         & \Phi_{l,c(m + 1) - 2} \circ \hdots \circ \Phi_{l,cm} \\
         & & \ddots
    \end{pmatrix}.
\end{equation}
With all algorithmic components defined, the two-level variant of \texttt{MGRIT} implemented in \texttt{XBraid} is given for a spatial discretization on $\Omega_l$ by \cref{alg:xbraid}. Here, the vectors on the coarse and fine scale are of different sizes and, of course,  have their corresponding lengths. 

\begin{algorithm}
    \caption{Two-Level cycle of \texttt{MGRIT} for fixed $\Omega_l$ implemented in \texttt{XBraid}}
    \label{alg:xbraid}
    \begin{algorithmic}[1]
	    \State Relax $B^{(N)}_l \vec{u}_l = \vec{g}^{(N)}_l$ using $F$-relaxation followed by $n_\text{relax} \geq 0$ applications of $CF$-relaxation
        \State Restrict fine grid solution and residual
        \begin{equation*}
		\vec{u}^{(N)}_{\text{coarse},l} \gets R^{(N)}_l \vec{u}_l, \quad \vec{r}^{(N)}_{\text{coarse},l} \gets R^{(N)}_l (\vec{g}^{(N)}_l - B^{(N)}_l \vec{u}^{(N)}_l)
        \end{equation*}
	    \State Solve $B^{(N)}_{\text{coarse},l} \vec{v}^{(N)}_{\text{coarse},l} = B^{(N)}_{\text{coarse},l} \vec{u}^{(N)}_{\text{coarse},l} + \vec{r}^{(N)}_{\text{coarse},l}$\label{alg:xbraid_solve_step}
	    \State Compute the coarse grid error $\vec{e}^{(N)}_{\text{coarse},l} \coloneq \vec{v}^{(N)}_{\text{coarse},l} - \vec{u}^{(N)}_{\text{coarse},l}$
	    \State Correct fine grid solution: $\vec{u}^{(N)}_l \gets \vec{u}^{(N)}_l + P^{(N)}_l \vec{e}^{(N)}_{\text{coarse}}$
    \end{algorithmic}
\end{algorithm}

Note that in practice the effective restriction is actually the transpose of the harmonic interpolation operator $P^{(N)}_l$, since restriction is always immediately preceded by $F$-relaxation.
Note also that \cref{alg:xbraid} does not exploit the linear nature of our model problem \cref{eq:parabolic_problem}. In \cref{alg:xbraid_solve_step} one could directly solve for $\vec{e}^{(N)}_{\text{coarse},l}$ due to the linearity of $B^{(N)}_{\text{coarse},l}$. This would make the restriction of $\vec{u}^{(N)}_l$ to $\vec{u}^{(N)}_{\text{coarse},l}$ unnecessary.
Nevertheless, we show the algorithm here as it is implemented in \text{XBraid} and presented in \cite{xbraid-package}. This more general FAS procedure allows the handling of nonlinear systems. However, this is not relevant for our work.

\subsubsection{Parallelization}
\label{sec:xbraid_parallelization}
\begin{figure}[htb]
	\centering
    \includegraphics[]{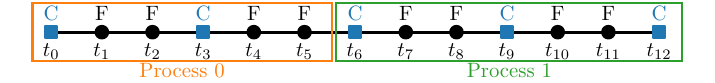}
    \caption{Distribution of $N = 12$ fine time steps onto $2$ processes for time coarsening factor $c=3$.}
    \label{fig:xbraid_distribution}
\end{figure}

\texttt{XBraid} allows the distribution of fine time steps across multiple processes in an \texttt{MPI} parallel fashion. The time steps are distributed in groups, where each $C$ time step $T_m = t_{cm}$ and all its subsequent $F$ time steps, i.e.\  $\{t_{cm}, t_{cm + 1}, \ldots, t_{(c + 1) m - 1}\}$, are on the same process. The last process also contains the last time step. An exemplary distribution is shown in \cref{fig:xbraid_distribution}.

\begin{figure}[htb]
	\centering
    \includegraphics[]{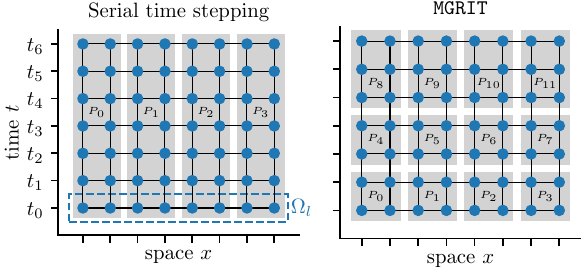}
    \caption{Space-time distribution onto processes. Figure reproduced and adapted from \cite{xbraid-package}.}
    \label{fig:xbraid_space_time_mpi}
\end{figure}

Due to its non-intrusive approach, \texttt{XBraid} is completely agnostic to the spatial decomposition. This means that the spatial dimension, i.e.\ the discretization $\Omega_l$ of $\Omega$, can be parallelized {\em independently} of the temporal parallelization. In traditional sequential time-stepping schemes, however, each process involved in the spatial parallelization is responsible for one part of the spatial discretization and its associated spatial data, such as matrix and solution entries, for all time steps.
These storage requirements can be reduced by noting that the spatial data is only needed for some of the time steps, depending on the chosen time propagator. For example, the backward Euler method requires only the spatial solution vector of the last time step and all spatial data for the current time step.
In contrast, the multigrid reduction-in-time approach requires the storage of all spatial data associated with the current process for all time steps, since each of the time steps can be solved repeatedly and without sequential order. Thus, the additional temporal parallelization level provided by \texttt{XBraid} ensures that each
process only needs to store its part of the time partition and thus only the corresponding spatial data associated with the process. This can be seen in \cref{fig:xbraid_space_time_mpi} for general time propagators, possibly depending on all previous time steps. 
When only using time coarsening, i.e.\ no space-time coarsening, the required memory per process for multigrid reduction-in-time increases by a factor of $O(\log_c N)$ with a time coarsening factor of $c > 1$ compared to sequential time stepping, see \cite{xbraid-package}.
In practice, \texttt{XBraid} supports memory usage reduction by storing the solution only at the C time steps and reconstructing solutions at F time steps by F relaxation whenever required, which is used in this work. We denote by $\hat{P}^t_l$ the temporal parallelization factor, i.e.\ the number of processes used for the temporal parallelization via \texttt{XBraid} for some spatial discretization $\Omega_l$. The total number of employed processes $P_l$ for the spatial and temporal parallelization of the problem associated with $\Omega_l$ is then given by 
\begin{equation}
\label{eq:subproblem_number_of_processes}
P_l = P^x_l\cdot \hat{P}^t_l,
\end{equation}
such that the overall number of processes is
\begin{equation}
\label{eq:total_number_of_processes}
P = \sum_l P_l.
\end{equation}

\section{A parallel-in-time combination method}
\label{sec:pitct}
In the following, we exploit the additional level of parallelization in the temporal component provided by the multigrid reduction-in-time scheme of \texttt{XBraid} in the context of the combination method and with our domain decomposition method for each of the subproblems arising there.
There are two prominent approaches: On the one hand, multigrid reduction-in-time can be used at the sparse grid level, resulting in a single global parallelization of time.
On the other hand, multigrid reduction-in-time can be applied within the combination method, i.e.\ it can be used to parallelize time for each subproblem of the combination method independently. In the following we will focus on the latter, since it allows, compared to the former, fine-grained control of the time step sizes for each subproblem mostly independently of each other with only loose coupling at the recombination step. This enables more temporal parallelism due to the fact that each subproblem can be parallelized independently. Additionally, the multigrid reduction-in-time on the subproblem level permits the use of a subproblem queue, i.e.\  between each of the recombination steps the subproblem solutions can be computed in any order, yielding independent compute tasks. This also allows the use of compute hardware with less compute capabilites, since the independent subproblem tasks can be managed by a job scheduling system such as SLURM (\cite{Jette2002}). In contrast to this, the naive global multigrid reduction-in-time approach produces one large compute task, however with the ability to recombine at any point in time.


\subsection{Multigrid reduction-in-time on the subproblems}\label{sec:parallel_in_time_subproblem_grid}
\begin{figure}[htb]
    \centering
    \includegraphics[]{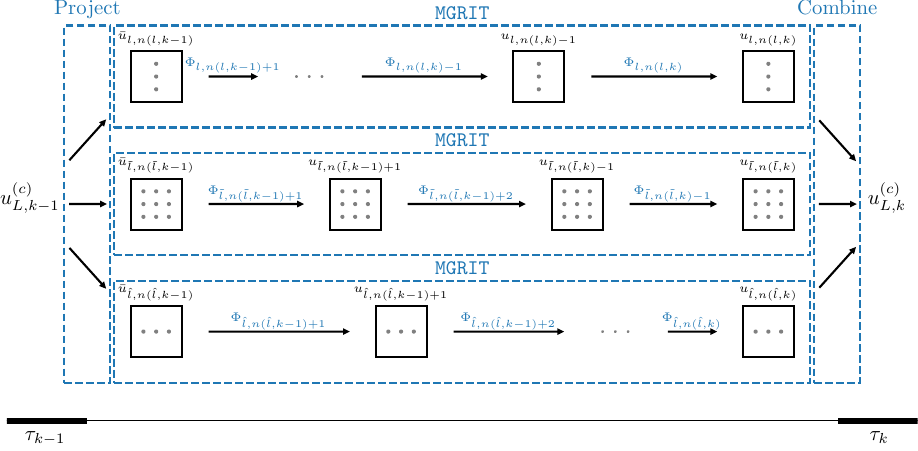}
    \caption{Example of one step of the combination method with multigrid reduction-in-time on the subproblems ($CTMGRIT^{loc}$).}
    \label{fig:ct_mgrit_subproblem_grid}
\end{figure}

In the following, we apply multigrid reduction-in-time to each subproblem in the sparse grid combination method individually.
For this, consider \cref{alg:combination_technique_general_parabolic}, where the time propagators $\{\Phi_{l, n(l, k - 1) + 1}, \ldots, \Phi_{l, n(l, k)}\}$ are used to sequentially propagate the solution from recombination step $\tau_{k-1}$ to $\tau_k$ for $1 \leq k \leq s$ and all subproblems $l$.
Since the time propagators $\Phi_{l, n}$ are independent of each other across subproblems, we can apply multigrid reduction-in-time to the time interval $(\tau_{k-1}, \tau_k]$ for each subproblem individually and simultaneously. This replaces the sequential nature of the solution scheme in time, i.e.\ the sequential application of $\{\Phi_{l, n(l, k - 1) + 1}, \ldots, \Phi_{l, n(l, k)}\}$, by an additional level of parallelism.
The resulting method is described in \cref{alg:ct_mgrit_subproblem_grid}, an example is given in \cref{fig:ct_mgrit_subproblem_grid}. Note also that, analogous to \cref{alg:combination_technique_general_parabolic}, the projection and combinations in \cref{alg:ct_mgrit_subproblem_project,alg:ct_mgrit_subproblem_combine} can be performed simultaneously as described in \cref{sec:hierarchisation}.
\begin{algorithm}
    \caption{$CTMGRIT^{loc}$: Combination method with multigrid reduction-in-time on the subproblems}
    \label{alg:ct_mgrit_subproblem_grid}
    \hspace*{\algorithmicindent} \textbf{Input} Spatial subproblem grids $\Omega_l$ for the combination method with sparse grid level $L$\\
    \hspace*{\algorithmicindent} \hphantom{\textbf{Input}} Time propagators $\Phi_{l, k}$ between recombination steps $\tau_k, 1 \leq k \leq s$\\
    \hspace*{\algorithmicindent} \hphantom{\textbf{Input}} Initial value $\bar{u}_0$\\
    \hspace*{\algorithmicindent} \textbf{Output} Approximate sparse grid solution $u^{(c)}_{L, s}$ at time step $T_s$
    \begin{algorithmic}[1]
            \State Set $u_{l, n(l, 0)}$ to the initial value $\bar{u}_0$ for all subproblems $l$.
	    \State Combine all $u_{l, n(l, 0)}$ using \cref{eq:combination_technique_solution} to obtain $u^{(c)}_{L, 0}$. \Comment{inter-subproblem comm.}
            \For{$1 \leq k \leq s$}
                \For{each subproblem $l$ simultaneously}
                    \State Project $u^{(c)}_{L, k - 1}$ onto $\Omega_l$ to obtain $\bar{u}_{l, n(l, k - 1)}$. \label{alg:ct_mgrit_subproblem_project}
                    \State Compute $u_{l, n(l, k)}$ from $\bar{u}_{l, n(l, k-1)}$ via \cref{alg:xbraid} 
		    \State $\quad \quad$ with the time propagators $\{\Phi_{l, n(l, k-1) + 1},
                \ldots, \Phi_{l, n(l, k)}\}$ \label{alg:ct_mgrit_subproblem_compute}
		\State Combine all $u_{l, n(l, k)}$ using \cref{eq:combination_technique_solution} to obtain $u^{(c)}_{L, k}$. \Comment{inter-subproblem comm.}\label{alg:ct_mgrit_subproblem_combine}
            \EndFor
        \EndFor
        \State\Return $u^{(c)}_{L, s}$
    \end{algorithmic}
\end{algorithm}

\cref{alg:ct_mgrit_subproblem_grid} adds a level of parallelism at the subproblem scale: For each subproblem, the problems associated with each time interval $(\tau_{k - 1}, \tau_k]$ are solved simultaneously and independently using the multigrid reduction-in-time scheme. 
The spatial problem can be parallelized again by the domain decomposition approach of \cref{sec:spatial_subproblem_solver}, resulting in a third level of parallelization.
We call this approach $CTMGRIT^{loc}$ and will use it throughout our numerical experiments in \cref{sec:numerical_experiments}.

Note that the choice of the number of recombination steps $s$ is highly problem dependent. In particular, to avoid subproblem solutions drifting apart for problems with large anisotropy, either frequent recombination of subproblem solutions, a significant reduction of the time step size, or both are required. However, due to the fact that an application of \texttt{MGRIT} is most impactful when many timesteps can be parallelized over, a tradeoff must be found. This tradeoff however must be found on a per-problem basis and is not the focus of this work. In the following we mitigate this problem by employing an initial level $L_0$ as introduced in \cref{sec:ct_fundamentals}, which reduces the anisotropy of the subproblems and allows for fewer recombination steps.

\begin{remark}
Note that applying MGRIT to the whole problem has the advantage that the subproblems can be recombined at any time. However, its major disadvantages are on the one hand the fixed
time steps for each subproblem, i.e. each subproblem has to obey the same time partition. This is particularly problematic due to the range of anisotropy of the subproblem grids, which then forces the time step size to be quite small for all subproblems and not only those highly anisotropic subproblems that require it. This is not the case for MGRIT on the subproblem level, i.e. $CTMGRIT^{loc}$, since there the time partitions of each subproblem only have to coincide at the recombination steps. On the other hand, applying MGRIT to the whole problem produces one large compute task, i.e. all subproblems have to be computed simultaneously. For $CTMGRIT^{loc}$ this is not the case, since the subproblems do not necessarily have to be computed simultaneously. Recombination could, for example, be done via checkpointing, i.e. writing solutions to file and running a recombination compute task. This allows for better load balancing via a job scheduling system as well as the usage of smaller compute systems.
\end{remark}

\section{Numerical experiments}
\label{sec:numerical_experiments}
Now we discuss the results of our numerical experiments for the parallel $CTMGRIT^{loc}$ approach. First, we give the parameters of our solver used in the experiments and describe the high-performance computer system we used. Then we show the parallelization properties of $CTMGRIT^{loc}$ by applying it to the heat equation,  the chemical master equation and some exemplary stochastic differential equations.

\subsection{Solver parameters and computer system}
\label{sec:solver_parameters}
Each of the components of $CTMGRIT^{loc}$ can be run with different parameter values and settings. Here, we do not intend to find an optimized parameter set\footnote{This is the subject of future research.} for any particular problem under consideration, but we rather settle on a parameter set that results in good performance for all problems considered. 
A spatio-temporal infinity norm with tolerance $10^{-8}$ is used, i.e.\ $CTMGRIT^{loc}$ terminates each \texttt{MGRIT} block whenever all spatial $l_2$ norms are less than $10^{-8}$. We employ the two-level version of \texttt{MGRIT} with time coarsening factor $c = 2$ within $CTMGRIT^{loc}$ as well as $FCF$-relaxation. The choice of time coarsening factor is motivated by the constraints on the maximal time step size $\Delta t$ imposed by the physics of each problem. Since a parameter study with respect to the optimization of runtime is not the focus of this work, the smallest possible $c$ was chosen to allow for higher dimensional problems on the available hardware. 
 
For problems with a symmetric operator as in \cref{sec:experiment_heat_equation}, the spatial solver of each subproblem itself is given by the conjugate gradient method preconditioned by our domain decomposition preconditioner using the balanced approach from \cref{eq:preconditioner_balanced} described in \cref{sec:spatial_subproblem_solver}. All spatial problems are solved up to a $l_2$ norm of the residual of $10^{-8}$.
The overlap parameter $\gamma$ is chosen as $\frac{1}{2}$, the coarse grid parameters $q_l$ are set to $q_l \equiv q \coloneq 2^{S - 4}$, where $2^S$ is the size of the spatial subdomains $\bar{\Omega}_{l, i}$ from \cref{eq:spatial_speedup_factor_new}.
All overlapping spatial subproblems and the spatial coarse grid problem are solved by a direct solver via the LU decomposition of \texttt{MUMPS} \cite{MUMPS:1,MUMPS:2} through the \texttt{PETSc} solver suite \cite{Dalcin2011,Balay2025}.
Furthermore, $S$ and thus $P^x_l$ from \cref{eq:spatial_speedup_factor_new} are chosen on a problem-by-problem basis.

The spatial solver for non-symmetric operators is given by the biconjugate gradient stabilized method, which is preconditioned by our domain decomposition method using the additive approach, where the LU decomposition is again used for the overlapping spatial subproblems and the spatial coarse grid problem. All other parameters are chosen as in the symmetric case.
The advantage of the additive approach in the non-symmetric setting is the reduction of memory by a factor of two as well as a reduction of communication, since the transposed operator does not have to be formed and stored.
In the symmetric setting, the transposed operator is simply given by the operator itself, so there is no additional cost. A comparison of the performance of the additive and the balanced approach for non-symmetric problems will be future work.

All computations have been run on the {\em Yuma} cluster of the Institut für Numerische Simulation of the University of Bonn. It consists of $45$ nodes, each equipped with two AMD EPYC 9435 32-Core processors and $768$ GiB of physical memory. The cluster employs an InfiniBand interconnect with $200$ Gb/s throughput. In the following numerical experiments we employ a one-to-one mapping of processes to cores. Therefore, each subproblem of $CTMGRIT^{loc}$, which corresponds to a single compute job on the cluster, can by parallelized by at most $45 \cdot 2 \cdot 32 = 2880$ processes. Hence, our choice of parallelization factors is constrained by the hardware limitations to $\hat{P}^t_l P^x_l \leq 2880$. Note however that this is only a limitation for the parallelization of each individual subproblem. Since the subproblems are independent of each other, the corresponding compute jobs can be processed in any order, for example via the SLURM job scheduling system \cite{Jette2002} employed on {\em Yuma}.

\subsection{The heat equation in higher dimensions}
\label{sec:experiment_heat_equation}
We first demonstrate the scalability properties of $CTMGRIT^{loc}$ described in \cref{sec:parallel_in_time_subproblem_grid}. To do this, we consider the standard heat equation
\begin{equation}
    \label{eq:heat_equation}
    \begin{aligned}
        \frac{\partial u}{\partial t} - \Delta u &= f\quad&&\text{in }\Omega\times(0, T],\\
        u &= 0\quad&&\text{on }\partial\Omega\times(0, T],\\
        u &= u^*\quad&&\text{in }\Omega\times\{0\},
    \end{aligned}
\end{equation}
where $\Omega=[0, 1]^d$ for $d=2,\ldots,6$, $T = 1$ and $u^*(x, t)=\sqrt{\norm{x}_2^2 + t^2}e^{-t}\prod_{i=1}^d\sin(\pi x_i)$. Here, $f$ is chosen so that $u^*$ is the exact solution.

\begin{figure}[htb]
	\centering
    \includegraphics[]{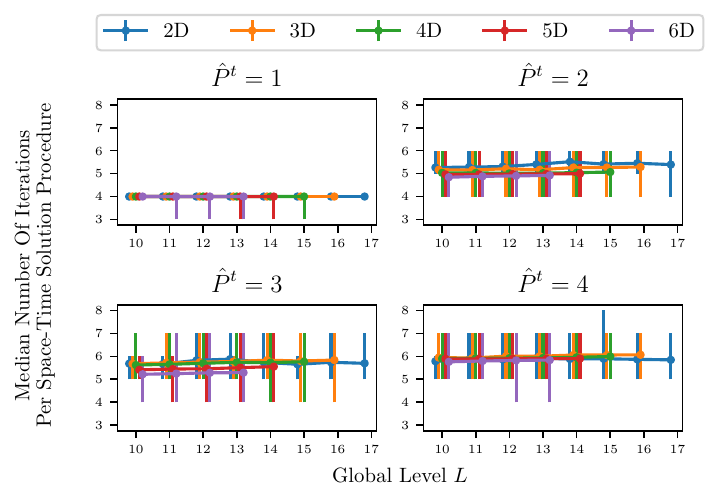}
\caption{Median number of iterations for all \texttt{MGRIT} solves for $CTMGRIT^{loc}$ for various choices of $d$ and $\hat{P}^t$. Note that the variation in the observed number of iterations (which is essentially bounded by a factor of $2$) is essentially due to a slight dependence of the MGRIT convergence rate on the anisotropy of the respective subproblem grid.}
\label{fig:heat_xbraid_iterations}
\end{figure}

In the notation of \cref{sec:general_algorithm}, we solve \cref{eq:heat_equation} by $CTMGRIT^{loc}$ on the sparse grid level $L$ using $s + 1$ recombination steps $\tau_k = k\frac{T}{s}, 0 \leq k \leq s,$ and a uniform time partition $t_{l, n} = t_n = n\frac{T}{N}, 0 \leq n \leq N_l = N$ in all subproblems $l$. Employing a subproblem-independent temporal parallelization factor $\hat{P}^t_l = \hat{P}^t$, see \cref{sec:xbraid_parallelization}, $N$ is chosen as $N \coloneq \hat{P}^t \cdot 10 \cdot s$, such that $10$ time steps are computed on each process between every two recombination steps for each subproblem.
The overlap factor $\gamma$ is set to $\gamma = 0.5$ and the number of processes for each subproblem $P^x_l$ is chosen according to \cref{eq:spatial_speedup_factor_new} with $S = 10$ so that each subdomain $\Omega_{l, i}$ of each subproblem $l$, and thus each processor, contains approximately $(1 + 2\gamma)2^S = 2\cdot2^{10} = 2^{11}$ grid points, see \cref{sec:spatial_subproblem_solver}.
The maximum sparse grid level $L_{\text{max}, d}$ in each dimension $d$ was set to $L_{\text{max}, d} = 19 - d$. Additionally, our choice of $N$ ensures that each processor is responsible for storing the spatial data associated with its grid nodes for $10$ time steps. Overall, this results in a balanced distribution of data across all processes involved. The total number of processes across all subproblems for these choices can be computed according to \cref{eq:spatial_speedup_factor_new,eq:subproblem_number_of_processes,eq:total_number_of_processes} and is given in \cref{tab:num_processes}.
\begin{table}[tb]
    \centering
    \caption{Total number of processes utilized by $CTMGRIT^{loc}$ for $L_{\text{max}, d} = 19 - d$, $S = 10$ and $\hat{P}^t_l = 1, \ldots, 4$. Recall that these process counts indicate the maximum number of processors that can be employed in parallel by solving all subproblems simulanteously, compare \cref{sec:solver_parameters}. }
    \label{tab:num_processes}
    \begin{NiceTabular}{c | c | c | c | c | c | c | c | c c }
        \diagbox{$d$}{$L$} & $10$ & $11$ & $12$ & $13$ & $14$ & $15$ & $16$ & $17$ & \\
        \cline{1-9}
        \phantom{a}2\phantom{a}& $27 $ & $56 $ & $120 $ & $260 $ & $564 $ & $1220 $ & $2628 $ & $5636$ &  \\
        3 & $247 $ & $547 $ & $1282 $ & $3024 $ & $7143 $ & $16743 $ & $38966 $ & & \\
        4 & $1249 $ & $3085 $ & $8017 $ & $20796 $ & $53576 $ & $136452 $ & & & \\
        5 & $4626 $ & $12842 $ & $36737 $ & $104522 $ & $292727 $ & & & & \\
        6 & $14152 $ & $43870 $ & $137608 $ & $425610 $ & & & & \phantom{$5636$} & \\
       \CodeAfter
        \begin{tikzpicture}
	       \draw[brace]
            ([xshift=3mm]2-9.north east) to node[auto = left] {$\cdot\hat{P}^t_l$} ([xshift=3mm]6-9.south east);
       \end{tikzpicture}
    \end{NiceTabular}
\end{table}
In the following we call the solution procedures between recombination steps for the subproblem $l$, i.e.\ \cref{alg:ct_mgrit_subproblem_compute}, the space-time solution procedure for the subproblem $l$. 
\cref{fig:heat_xbraid_iterations} shows the median number of iterations for the space-time solution procedure across all subproblems.
We can clearly observe an optimal scaling behavior of $CTMGRIT^{loc}$, independent of dimension $d$, global level $L$, subproblems $l$ and therefore anisotropic subproblem grids.

\begin{figure}[htb]
	\centering
    \includegraphics[]{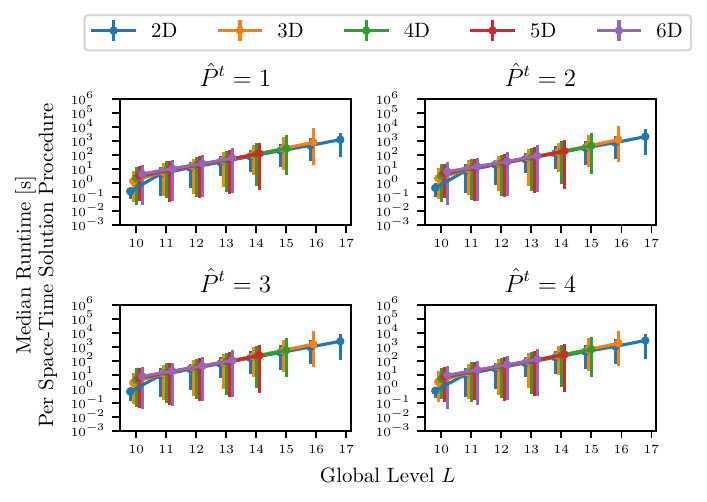}
\caption{Median runtime in seconds for all \texttt{MGRIT} solves for $CTMGRIT^{loc}$ for various choices of $d$, $L$ and $\hat{P}^t$. The runtimes are given on a logarithmic scale so that small variations in the observed runtimes are barely visible. Note also that the total number of time steps $N_l=N$ employed is independent of $l$ and is linearly dependent on the temporal parallelization factor $\hat{P}_t$ so that similar runtimes for different values of $\hat{P}_t$ indicate an optimal weak scaling behavior. Moreover, we can observe only a very minor dependence of the runtimes on the spatial dimension $d$.}
\label{fig:heat_xbraid_times}
\end{figure}

\cref{fig:heat_xbraid_times} shows the median runtime for the space-time solution procedure over all subproblems. In view of the scaling experiments for the spatial solver within the combination method performed in \cite{Griebel2023}, we can now observe that the scaling behavior of $CTMGRIT^{loc}$ is controlled by the scaling behavior of the spatial solver, since it shows qualitatively very similar runtimes for each subproblem. The significantly faster space-time solution procedure for $d=2$ and $L=10$ is due to the fact that, for this combination of parameters, each spatial subproblem is distributed over at most two processes. In this particular case, since the spatial subproblems are sufficiently small and the overlap factor was chosen to be $\gamma = 0.5$, the spatial domain decomposition solver is essentially a direct solver, since the full problem is available for each process. This allows the spatial subproblem to be solved in very few iterations.

\begin{figure}[htb]
	\centering
    \includegraphics[]{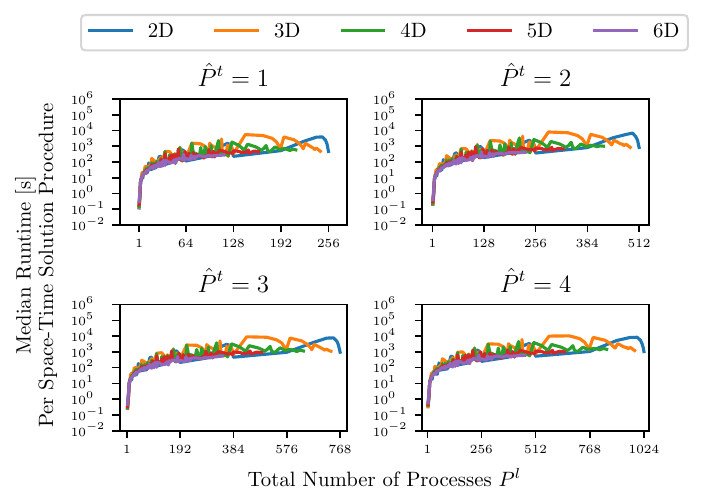}
\caption{Median runtime in seconds for all \texttt{MGRIT} solves for $CTMGRIT^{loc}$ compared to the total number of processes. The runtimes are given on a logarithmic scale.}
\label{fig:heat_xbraid_scaling}
\end{figure}

The data of \cref{fig:heat_xbraid_times} can also be interpreted as a weak-scaling study. This is due to our specific choice of spatial and temporal parallelization factors via \cref{eq:spatial_speedup_factor_new,eq:subproblem_number_of_processes} such that the amount of work per process is approximately constant. \cref{fig:heat_xbraid_scaling} shows the median runtime against the total number of processes. We observe near optimal weak scaling after a preasymptotic phase, which are the same results as for our spatial domain decomposition solver presented in \cite{Griebel2023}. We can therefore deduce that $CTMGRIT^{loc}$ inherits its scaling properties from the spatial solver.

\subsection{An approximation to chemical and biological reaction networks}
Our next problem is given by the so-called chemical master equation (CME). It models the evolution of the probability distribution of the configuration of a reaction network, where in particular a configuration describes the number of molecules of different species, such as chemical reactants, undergoing reactions. Examples range from classical stochastic chemical kinetics or gene expression in a population of cells to epidemiological models. The chemical master equation is often solved by a reduction of its very large state space. Prominent examples are the Gillespie algorithm \cite{Gillespie1977} (also known as SSA) based on a stochastic reduction of the state space or projection algorithms such as the Krylov-based finite state projection algorithm of \cite{Burrage2006}. In the following we will demonstrate the application of $CTMGRIT^{loc}$ to the chemical master equation, which circumvents a reduction of the state space via its inherent sparse grid approximation properties.

We follow the presentation in \cite{Lunz2021}. Consider a biochemical system with $d$ species $\{S_1, \ldots, S_d\}$.
Let $\nu_r\in\mathbb{Z}^d$ and $\alpha_r \geq 0$ denote the stoichiometric vector and the propensity for each reaction $1\leq r \leq m$ between these species.
In particular, the reaction $r$ can be written as
\begin{equation}
    a_1S_1 + \ldots + a_dS_d \xrightarrow{\alpha_r(x)} b_1S_1 + \ldots + b_dS_d,
\end{equation}
i.e.\ the current state $x=(x_1, \ldots, x_d)\in\mathbb{N}_{\geq 0}^d$ is changed by the reaction $r$ to the state $x + \nu_r$ with the probability per unit time given by $\alpha_r(x)$, where $\nu_r = (b_1 - a_1, \ldots, b_d - a_d)$.
It is well known that this process can be modeled by a Markov chain on the integer lattice whose probability density function satisfies the chemical master equation
\begin{equation}
    \label{eq:cme}
    \frac{\partial u(x, t)}{\partial t} = \sum_{r=1}^m \alpha_r (x - \nu_r) u(x - \nu_r, t) - \alpha_r(x)u(x, t),
\end{equation}
for some initial condition $u(x, 0) = u_0(x)$.
Via an embedding of the integer lattice in $\Omega \subset \mathbb{R}^d$ and the Kramers-Moyal expansion, a Fokker-Planck-type approximation to the chemical master equation \cref{eq:cme} can be written as
\begin{equation}
    \label{eq:cme_fp}
    \frac{\partial u(x, t)}{\partial t} = \sum_{r=1}^m -\nu_r^T \nabla\left[\alpha_r(x)u(x, t)\right]
            + \frac{1}{2}\nu_r^T (\nabla \otimes \nabla\left[\alpha_r(x)u(x, t)\right])\nu_r,
\end{equation}
which matches our parabolic problem \cref{eq:parabolic_problem}.
By making the computational domain $\Omega$ sufficiently large such that the mass is far from the boundary, we can assume homogeneous Dirichlet boundary conditions on $\partial\Omega$.

\subsubsection{A genetic toggle switch in two dimensions}
\begin{figure}[htb]
    \centering
    \begin{subfigure}[c]{0.48\textwidth}
        \centering
        \includegraphics[width=\textwidth]{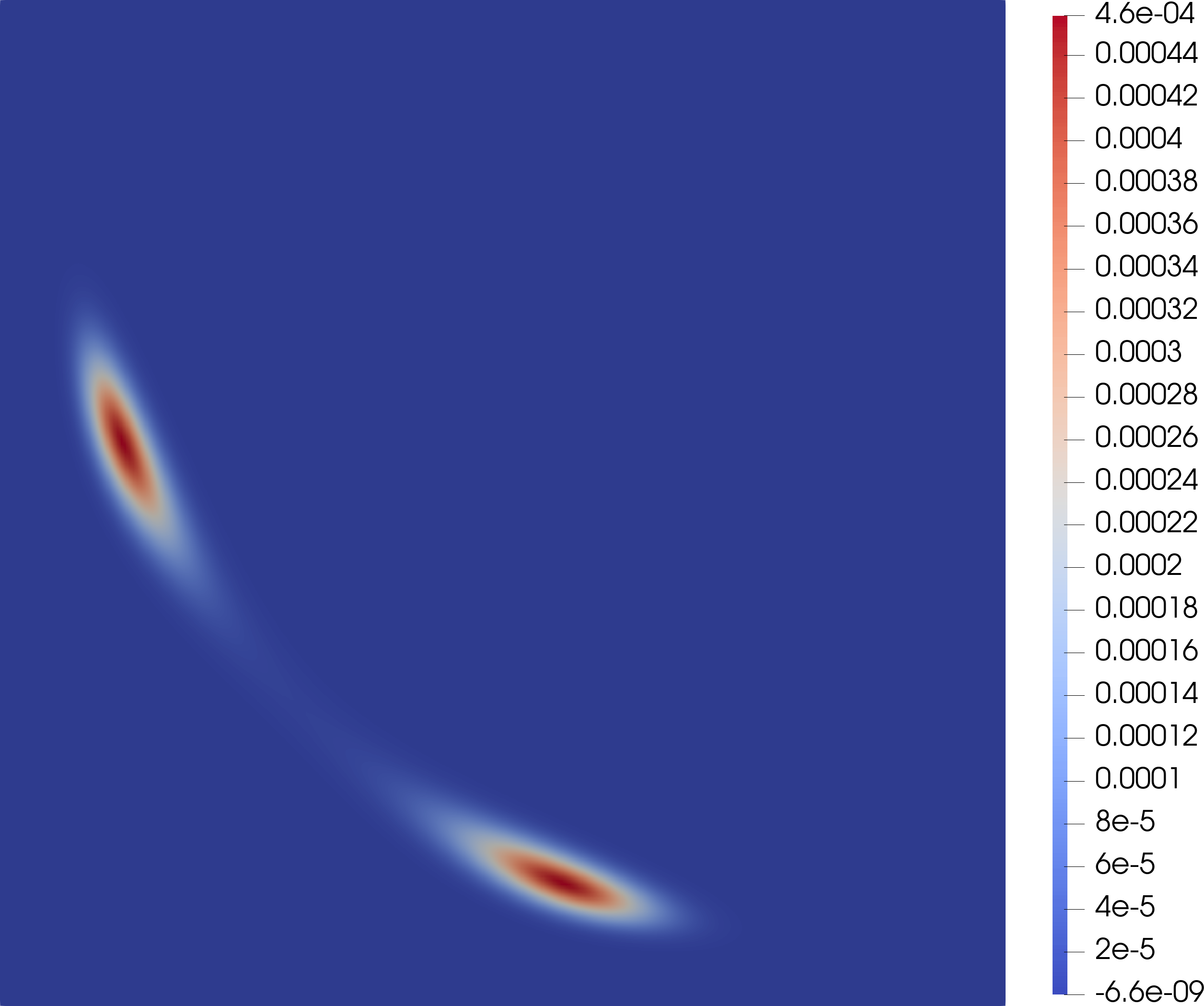}
    \end{subfigure}%
    \hfill%
    \begin{subfigure}[c]{0.48\textwidth}
        \centering
        \includegraphics[width=\textwidth]{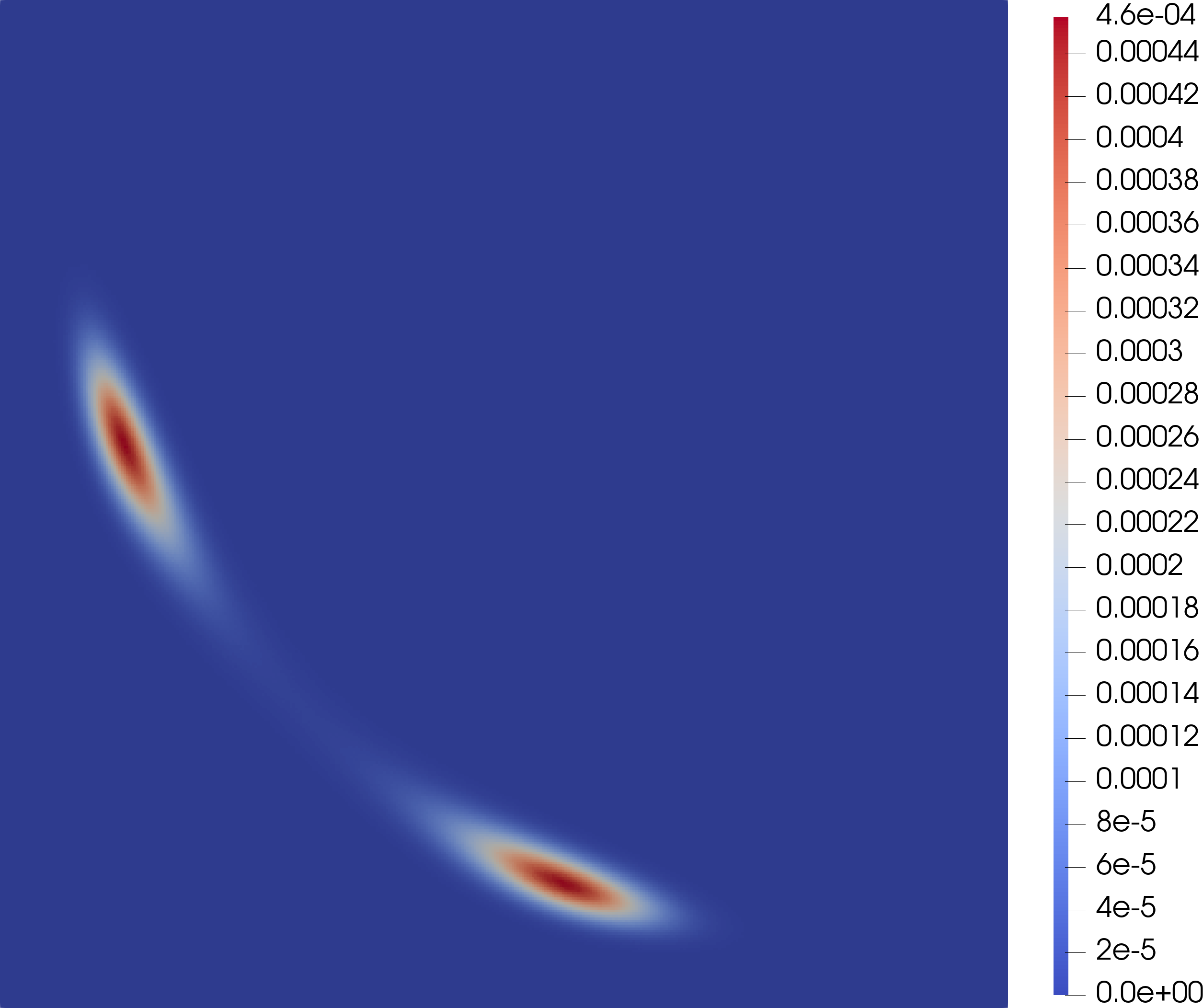}
    \end{subfigure}
    \begin{subfigure}[c]{0.48\textwidth}
        \includegraphics[width=0.82\textwidth]{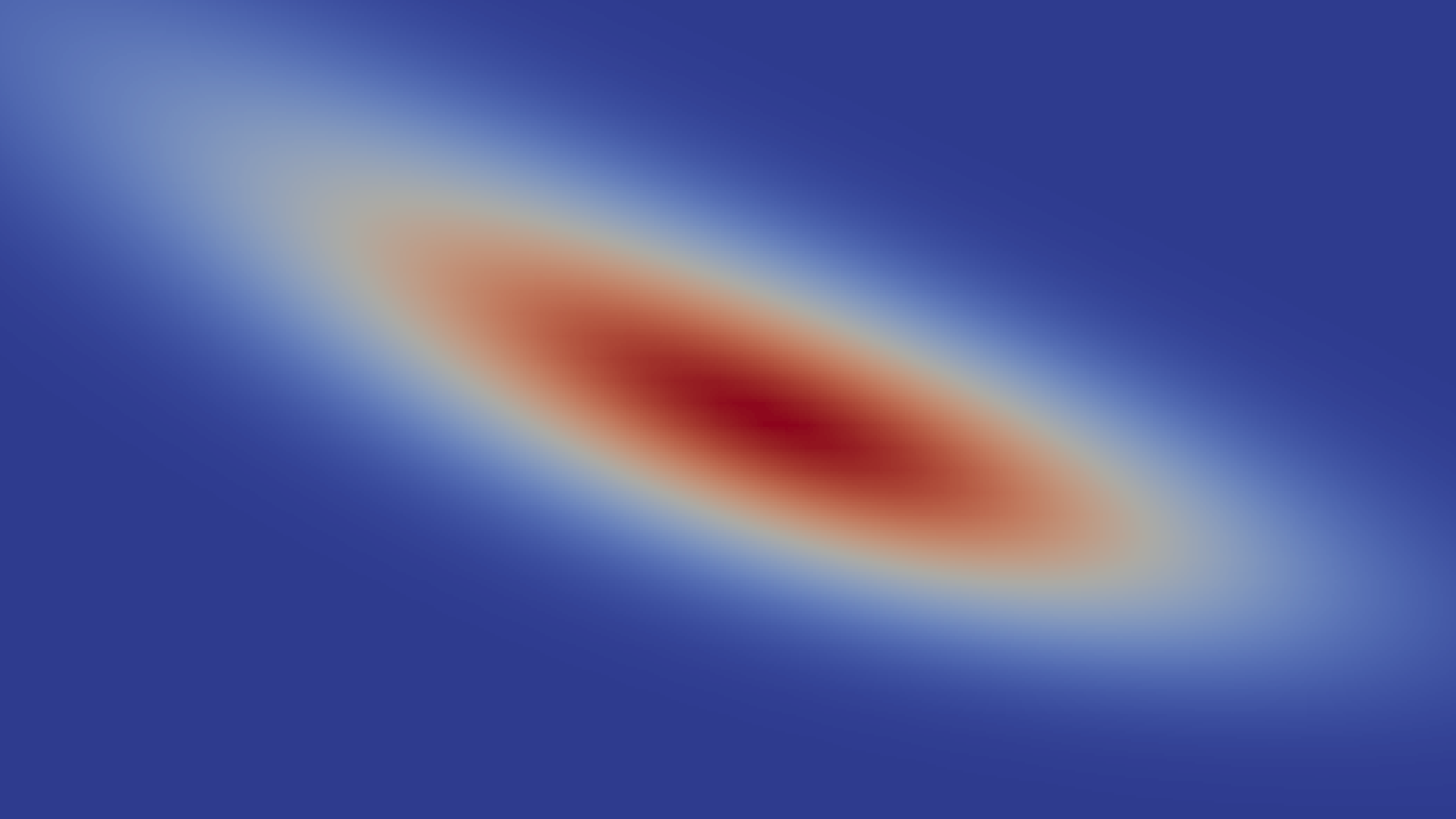}
    \end{subfigure}%
    \hfill%
    \begin{subfigure}[c]{0.48\textwidth}
        \includegraphics[width=0.82\textwidth]{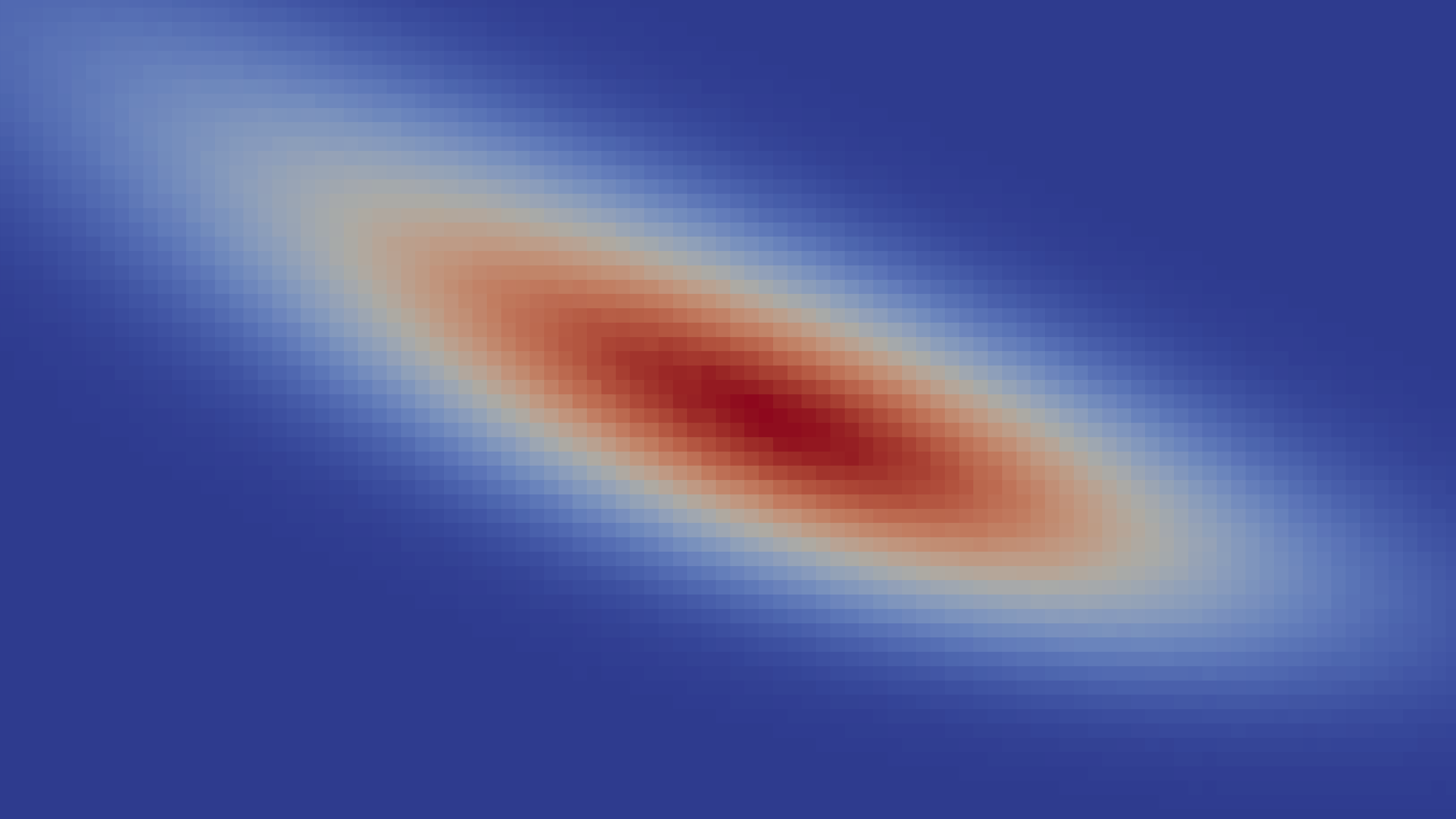}
    \end{subfigure}
    \caption{Approximate solution (top) for the two-dimensional toggle switch computed via $CTMGRIT^{loc}$ (left) and SSA (right) at time $t = 10^5s$ with closeups (bottom).}
	\label{fig:toggle_switch_2d_solutions}
\end{figure}

As an example we consider a genetic toggle switch. It models two competing repressors $A$ and $B$, which are transcribed by two constitutive promoters and can each inhibit the production of the other repressor. An example for this is the genetic toggle switch in Escherichia coli in \cite{Gardner2000}. The model exhibits a bistable stationary distribution and has been used extensively in the literature to test numerical schemes for the chemical master equation, see \cite{Kryven2015,Deuflhard2008,Sjoeberg2007} among others. Due to its multi-stability, the convergence of traditional stochastical schemes such as the Gillespie algorithm \cite{Gillespie1977} can be slow to approximate stationary distributions \cite{Kryven2015}. In this context, deterministic numerical schemes based on the chemical master equation such as our proposed method are clearly advantageous. The four reactions and their propensities used in the following are given in \cref{tab:toggle_switch_2d}, where the species $\{S_1, S_2\}$ and the current state $(x_1, x_2)$ have been relabeled to  $\{A, B\}$ and $([A], [B])$ for better readability. Following \cite{Kryven2015} we use a computational domain of $[0, 399]^2$ and let the initial condition $u_0$ be given by the probability density function of a two-dimensional normal distribution with mean $M = \begin{bmatrix} 133 & 133\end{bmatrix}^T$ and variance $C = \diag(133)$. We run $CTMGRIT^{loc}$ on level $L = 13$ with initial level $L_0 = 6$, $S = 10$ and two recombination steps, i.e.\ $s = 1$. We set $\hat{P}^t = 100$ and compute $N = 100000$ time steps, such that each process handles $1000$ time steps. In total $CTMGRIT^{loc}$ consists of five independent spatial-temporal subproblems, three of which utilize $1600$ processes, the other $800$ processes, for a total number of $6400$ processes.
\cref{fig:toggle_switch_2d_solutions} shows the solution at time $t=10^5s$ of $CTMGRIT^{loc}$ as well as the probability density function extracted from $10^8$ trajectories generated by \textit{SSA}. We observe that $CTMGRIT^{loc}$ produces a smooth probability density, whereas the probability density extracted from the trajectories of \textit{SSA} has jumps, since \textit{SSA} is an integer-valued scheme. These jumps can only be smoothed by significantly increasing the number of sampled trajectories. The subproblems of $CTMGRIT^{loc}$ had a median runtime of $2959.34$ seconds, with a minimum of $2091.67$ seconds and maximum of $3682.76$ seconds. Each trajectory of \textit{SSA} was computed using the Python bindings of the Rust crate \texttt{rebop} of \cite{Andreani2024}. Each trajectory was simulated with a median runtime of $0.0045836$ seconds. Thus, for $10^8$ trajectories, with the same number of processes as for $CTMGRIT^{loc}$, the total runtime of the \textit{SSA} approach is $71.62$ seconds. However, due to the vastly different algorithmic nature of both approaches, and in particular approximation and convergence properties in more complicated scenarios, a comparison of runtimes is inconclusive.

\begin{table}
    \centering
    \caption{The two-dimensional genetic toggle switch model as described.}
    \label{tab:toggle_switch_2d}
    \begin{tabular}{r c l | r c l | l c}
        \multicolumn{3}{c|}{reaction} & \multicolumn{3}{|c|}{propensity} & \multicolumn{2}{c|}{rate}\\
        \hline
        A & $\to$ & 2A & $\alpha_1$ & $=$ & $c_1/(c_2 + [B]^\beta)$ & $c_1 = 3 \cdot 10^3, c_2 = 1.1 \cdot 10^4, \beta = 2$ & \\
        2A & $\to$ & A & $\alpha_2$ & $=$ & $c_3[A]$ & $c_3 = 10^{-3}$ & \\
        B & $\to$ & 2B & $\alpha_3$ & $=$ & $c_4/(c_5 + [A]^\gamma)$ & $c_4 = 3 \cdot 10^3, c_5 = 1.1 \cdot 10^4, \gamma = 2$ & \\
        2B & $\to$ & B & $\alpha_4$ & $=$ & $c_6[B]$ & $c_6 = 10^{-3}$ &
    \end{tabular}
\end{table}

\subsubsection{A genetic toggle switch in three dimensions}
\begin{figure}[htb]
    \centering
    \includegraphics[width=0.49\textwidth]{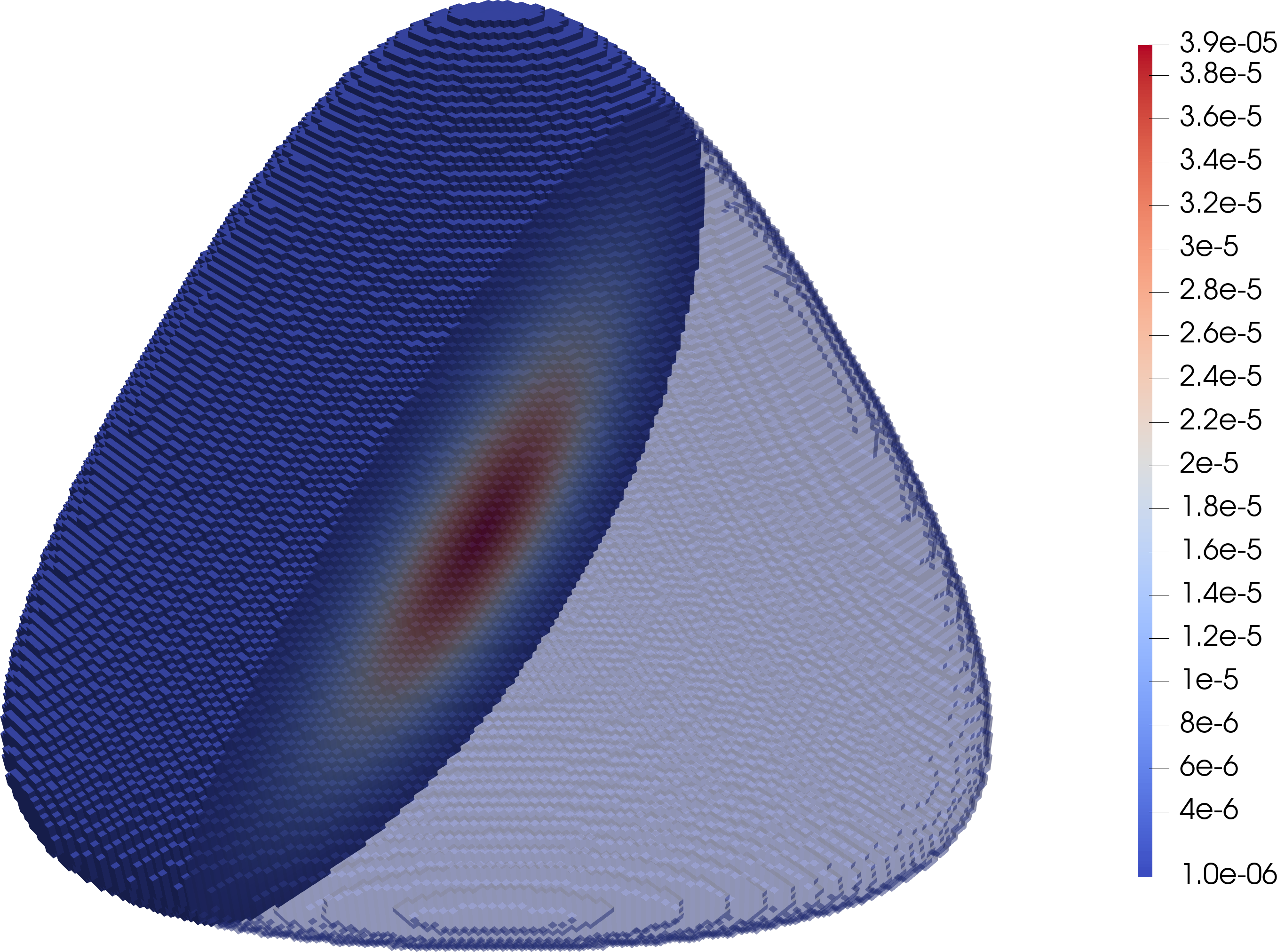}
    \includegraphics[width=0.49\textwidth]{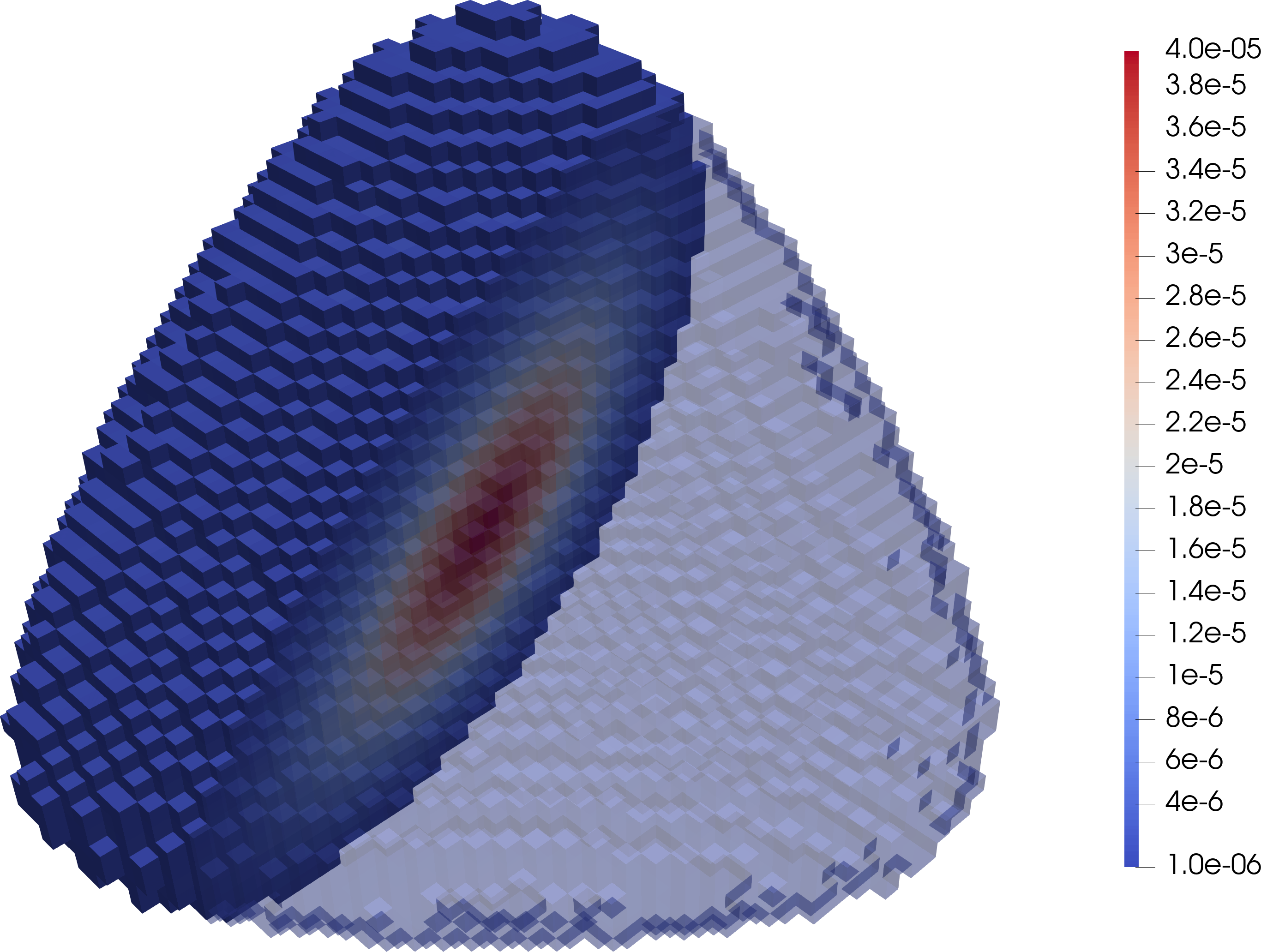}
    \caption{Approximate solution for the three-dimensional toggle switch computed via $CTMGRIT^{loc}$ (left) and SSA (right) at time $t = 10^5s$. The solution for $CTMGRIT^{loc}$ is sampled on a $300^3$ full grid.}
	\label{fig:toggle_switch_3d_solutions}
\end{figure}

\cite{Kryven2015} proposed an extension of the bistable toggle switch to a tristable toggle switch in three dimensions. \cref{tab:toggle_switch_3d} shows the reactions, where we have relabeled the species as $A, B$ and $C$. The computational domain is given by $[0, 199]^3$ and we use as the initial condition, analogously to the two dimensional problem, a three-dimensional Gaussian distribution centered at $\begin{bmatrix} 133, 133, 133\end{bmatrix}^T$ with variance $\diag(133)$. $CTMGRIT^{loc}$ is run on level $L = 20$ with initial level $L_0 = 7$, $S = 15$ and two recombination steps, i.e.\ $s = 1$. We set $\hat{P}^t = 10$ and compute $N = 10000$ time steps, such that each process handles $1000$ time steps. In total $CTMGRIT^{loc}$ consists of four independent spatial-temporal subproblems, three of which utilize $1260$ processes, the other $630$ processes, for a total number of $4410$ processes. All other parameters are the same as in the two-dimensional case. \cref{fig:toggle_switch_3d_solutions} shows the solutions at $t=10^5s$ of both $CTMGRIT^{loc}$ and \textit{SSA} with $10^8$ trajectories. We notice that $CTMGRIT^{loc}$ again produces a significantly more resolved probability density function. Note that the solution for $CTMGRIT^{loc}$ is sampled on a $300^3$ full grid due to the fact that it is not feasible to store the equivalent full grid on level $20$ in memory due to its size. The median runtime of each $CTMGRIT^{loc}$ subproblem is $3198.63$ seconds, whereas for the same number of processes the $10^8$ \textit{SSA} trajectories can be solved in $149.93$ seconds. However, analogously to the two dimensional scenario, the runtimes of both approaches are difficult to compare due to their different nature.

\begin{table}
    \centering
    \caption{The three-dimensional genetic toggle switch model as described.}
    \label{tab:toggle_switch_3d}
    \begin{tabular}{r c l | r c l | l c}
        \multicolumn{3}{c|}{reaction} & \multicolumn{3}{|c|}{propensity} & \multicolumn{2}{c|}{rate}\\
        \hline
        A & $\to$ & 2A & $\alpha_1$ & $=$ & $c_1/(c_2 + ([B] + [C])^\beta)$ & $c_1 = 3 \cdot 10^3, c_2 = 1.1 \cdot 10^4, \beta = 2$ & \\
        2A & $\to$ & A & $\alpha_2$ & $=$ & $c_3[A]$ & $c_3 = 10^{-3}$ & \\
        B & $\to$ & 2B & $\alpha_3$ & $=$ & $c_4/(c_5 + ([A] + [C])^\gamma)$ & $c_4 = 3 \cdot 10^3, c_5 = 1.1 \cdot 10^4, \gamma = 2$ & \\
        2B & $\to$ & B & $\alpha_4$ & $=$ & $c_6[B]$ & $c_6 = 10^{-3}$ & \\
        C & $\to$ & 2C & $\alpha_5$ & $=$ & $c_7/(c_8 + ([A] + [B])^\xi)$ & $c_7 = 3 \cdot 10^3, c_8 = 1.1 \cdot 10^4, \xi = 2$ & \\
        2C & $\to$ & C & $\alpha_6$ & $=$ & $c_9[C]$ & $c_9 = 10^{-3}$ &
    \end{tabular}
\end{table}

\subsection{Stochastic differential equations}
\label{sec:sde}
Next, we consider stochastic differential equations of the form
\begin{equation}
    \label{eq:sde}
    dX_t = \theta X_t \,dt + \sigma \,dW_t
\end{equation}
where $X_t$ is a $\mathbb{R}^d$-valued stochastic process, $\theta\in\mathbb{R}^{d \times d}$, $\sigma\in\mathbb{R}^{d \times m}$, and $W_t$ is an $m$-dimensional Gaussian white noise process with zero mean and autocorrelation $\mathbb{E}[W_{t + \tau} \overline{W_t}] = 2 D \delta(\tau)$ with $D \in \mathbb{R}^{m \times m}$.

The probability density function $u(x, t)$ of the Markov process $X_t$ can be computed using the Fokker-Planck equation in $\mathbb{R}^d$
\begin{equation}
    \label{eq:sde_fokker_planck}
    \frac{\partial u}{\partial t} = - \sum_{i=1}^d \frac{\partial}{\partial x_i} \left[(\theta x)_i u\right] + \frac{1}{2}\sum_{i=1}^d\sum_{j=1}^d \frac{\partial^2}{\partial x_i \partial x_j}\left[ H_{ij}u \right],
\end{equation}
with $H = 2 \sigma D \sigma^T$, boundary conditions $u\to 0$ as $\norm{x}\to\infty$ and a suitable initial condition $p(x, 0) = p_0(x)$.
The formal solution of \cref{eq:sde} is well known. It is given by
\begin{equation}
    \label{eq:sde_analytical_solution}
    X_t = e^{\theta t} X_0 + \int_0^t e^{\theta (t - s)}\sigma\,dW_{s}.
\end{equation}
and has a $d$-dimensional normal distribution $\mathcal{N}(e^{\sigma t} X_0, \mathrm{cov}(X_0, X_0) + \int_0^t e^{\theta s} H e^{\theta^T s} \,ds)$.

\subsubsection{A two-dimensional linear oscillator}
\label{sec:linear_oscillator_2d}
A first example for \cref{eq:sde_fokker_planck} is taken from \cite{Spencer1993}. It models the response of a two-dimensional linear oscillator  subject to additive Gaussian white noise. We consider, in the notation of \cref{eq:sde} and \cref{eq:sde_fokker_planck}, the system given by
\begin{equation}
\theta \coloneq \begin{bmatrix}
        0 & 1\\
        - \omega_0^2 & -2 \xi \omega_0
        \end{bmatrix}, \quad
\sigma \coloneq \begin{bmatrix}
                0 \\
                1
            \end{bmatrix}, \quad
D = \begin{bmatrix}0.1\end{bmatrix},
\end{equation}
where $\xi = 0.05$ and $\omega_0 = 1$. The system describes a mass-spring-damper model with a one-dimensional mass connected to a fixed point via a spring and a damper. An excitation $W_t$ is applied to the mass. The spatial dimensions of the system correspond to the position ($x_1$) and velocity ($x_2$) of the mass.
We choose a two-dimensional normal distribution with mean $M = \begin{bmatrix}5 & 5\end{bmatrix}^T$ and covariance $C = \diag(\frac{1}{9})$ as initial condition.
Using \cref{eq:sde_analytical_solution} and the fact that normal distributions with normal mean are normal, the analytical solution is given by a Gaussian process with mean $e^{\theta t} M$ and covariance $e^{\theta t} C e^{\theta^T t} + \int_0^t e^{\theta s} H e^{\theta^T s} \,ds$.

Since \cref{eq:sde_fokker_planck} is a simple two-dimensional problem, it can be solved with sufficient accuracy using a classical full-grid approach using a second order finite difference discretization. The discretization and parameters of the full-grid approach are chosen analogously to those of the subproblems of the combination method. The full-grid resolution is chosen such that the number of degrees of freedom is approximately the same as the sum of the degrees of freedom of all subproblems in the combination method.
This serves as a means to validate $CTMGRIT^{loc}$ not only against an analytical solution \cref{eq:sde_analytical_solution}, but also against a reference solution derived from a classical scheme.
In all cases we consider the computational domain $[-10, 10]^2$. The solution of the combination method has been computed on the sparse grid level $L=14$ with initial level $L_0=6$, i.e.\ $l\geq 6$ component-wise for all subproblems $l$. The full grid schemes employ a grid of size $421 \times 421$ such that the number of degrees of freedom approximately equals the total number of degrees of freedom across all subproblems of $CTMGRIT^{loc}$. 
We run the simulation until time $T_\text{end} = 100 \approx 16\frac{2\pi}{\omega_0}$, which corresponds to about $16$ rotations around the origin in phase space. We employ the backward Euler method as the time propagator $\Phi_{l, n}$ on the fine time partition with a uniform time step size $\Delta t =  0.005$. The time coarsening factor is chosen as $c = 2$ such that $\Delta T = 0.01$.
We compute the solution using the full grid with sequential time stepping (\textit{FGS}), the full grid with \texttt{XBraid} for time stepping (\textit{FGX}), the time-sequential combination technique (\textit{CT}) and $CTMGRIT^{loc}$ with two recombination steps, i.e.\ $s=1$.
For all three schemes we use our spatial domain decomposition solver with the parameters described in \cref{sec:solver_parameters}. We set $S=10$ and $q=6$ such that each process stores about $2048$ fine and $64$ coarse spatial degrees of freedom.


\Cref{tab:linear_oscillator_2d} shows the runtimes, total number of processes and maximum error over time on the spatial diagonal for each of the schemes. Note that the runtime is averaged over five repetitions. We can see that all schemes produce comparable solutions, where our $CTMGRIT^{loc}$ is both the most accurate and fastest. In fact, $CTMGRIT^{loc}$ significantly outperforms a naive application of \texttt{MGRIT} via \textit{FGX} for the same number of total processes used. Note that even though \textit{FGS} utilizes no temporal speedup, i.e.\ the number of processes is reduced by a factor of $10$, it only has a slightly slower runtime compared to $CTMGRIT^{loc}$. This is due to the fact that the size of the problem at hand is near the crossover point where the performance benefits of using a temporal parallelization via multigrid in time offsets the incurred overheads of such an approach. Additionally, there is an optimal number of spatial processes for the full grid scheme \textit{FGS} beyond which the performance starts to deteriorate again due to communication overheads. This can also be seen by running $\textit{FGS}$ with $1760$ spatial processes, the same number as total processes for $CTMGRIT^{loc}$, where we observed an average runtime of $61890$ seconds for \textit{FGS}. Therefore $CTMGRIT^{loc}$, or more generally multigrid in time schemes, allow to increase the performance beyond what is possible with a purely spatial parallelization scheme such as \textit{FGS}. 

The performance of $CTMGRIT^{loc}$ for a heuristically optimized set of parameters is shown in \cref{tab:linear_oscillator_2d_opt}. Due to hardware constraints and the large parameter space, we consider the reduced problem with $T_\text{end} = 1$ and $\delta t = \frac{1}{3200}$, such that there are $3201$ time steps. Only the two solvers based on the combination technique were considered due to their superior performance compared to the full grid variants. In this scenario we choose the residual mode of \texttt{XBraid} with a space-time residual tolerance of $1e-4$, such that our spatial solver is used as a preconditioner for the space-time system. We additionally employ two temporal levels with a coarsening factor of $c = 12$, choose $F$-relaxation and only apply two iterations of the spatial solver. All other parameters have not been changed. We can observe a substantial improvement in the performance of $CTMGRIT^{loc}$ compared to the  version, . An additional improvement not considered here is a subproblem specific choice of parameters for \texttt{XBraid} to mitigate the observed variation in the runtimes per subproblem. However, a full optimization of the parameters of each of the schemes to obtain the optimal performance is future work. We refer to, amongst others, \cite{Falgout2014} for a parameter study in the context of multigrid in time.

\begin{table}[htb]
    \centering
    \caption{Average, minimum and maximum runtimes per problem, total number of processes used and maximum error over time in the spatial origin for the two-dimensional linear oscillator.}
    \label{tab:linear_oscillator_2d}
    \begin{tabular}{c | l | l | c }
	    & run time [s] 	&  number of processes 	&  max. error in origin  \\
            \hline
            \textit{FGS} & $11595$ ($8656$/$14599$) & $176$ & 0.00098 \\
            \textit{FGX} & $66449$ ($66410$/$66489$) & $1760 = 10 \cdot 176$ & 0.00096 \\
            $CTMGRIT^{loc}$ & $10385$ ($6362$/$14744$) & $1760 = 10 \cdot (4 \cdot 32 + 3 \cdot 16)$ &  0.00085\\
            \textit{CT} & $1326$ ($721$/$1673$) & $176 = 4 \cdot 32 + 3 \cdot 16$ & 0.00089
    \end{tabular}
\end{table}

\begin{table}[htb]
    \centering
    \caption{Average, minimum and maximum runtimes per problem, total number of processes used and maximum error over time on the spatial diagonal for the two-dimensional linear oscillator with improved $CTMGRIT^{loc}$ parameters.}
    \label{tab:linear_oscillator_2d_opt}
    \begin{tabular}{c | l | l | c }
	    & run time [s] 	&  number of processes 	&  max. error on diagonal \\
            \hline
            $CTMGRIT^{loc}$ & $107$ ($62$/$156$) & $1760 = 10 \cdot (4 \cdot 32 + 3 \cdot 16)$ &  0.0028\\
            \textit{CT} & $162$ ($108$/$213$) & $176 = 4 \cdot 32 + 3 \cdot 16$ & 0.008
    \end{tabular}
\end{table}

\subsubsection{A four-dimensional linear system}
\begin{figure}[htb]
    \centering
    \includegraphics[]{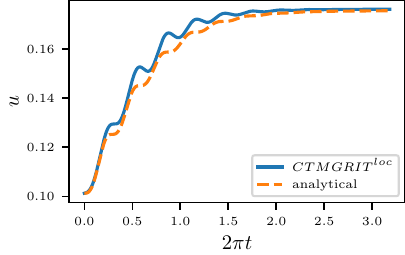}
    \includegraphics[]{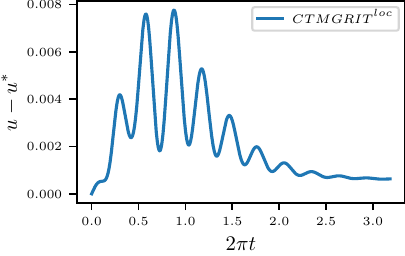}
    \caption{Approximate solution (left) and error (right) at the spatial origin over time of the four-dimensional linear oscillator for $CTMGRIT^{loc}$.}
	\label{fig:linear_oscillator_4d_origin}
\end{figure}

The second example of a stochastic differential equation \cref{eq:sde} is taken from \cite{Wojtkiewicz2001}. It is given by
\begin{equation}
\theta = \begin{bmatrix}
    0 & 1 & 0  & 0 \\
    -(k_1 + k_2) & -c1 & k2 & 0 \\
    0 & 0 & 0 & 1 \\
    k_2 & 0 & -(k_2 + k_3) & -c_2
\end{bmatrix},
\quad \sigma = \begin{bmatrix}
    0 & 0 \\
    1 & 0 \\
    0 & 0 \\
    0 & 1
\end{bmatrix},
\quad D = \begin{bmatrix}
    2D_1 & 0 \\
    0 & 2D_2
\end{bmatrix},
\end{equation}
where $k_1=k_2=k_3=1, c_1=c_2=0.4$ and $D_1=D_2=0.2$. The system describes a mass-spring-damper model with two one-dimensional masses. Each mass is connected to a separate fixed point via a spring and a damper. Additionally, the masses are connected via a third spring. An excitation given by $W_t$ is applied to both masses. The four dimensions of the system correspond to the positions ($x_1, x_3$) and velocities ($x_2, x_4$) of the masses.
The chosen initial condition is a zero mean Gaussian distribution with variance $C = \diag(0.5)$.
Analogous to \cref{sec:linear_oscillator_2d}, the analytical solution can be derived from the initial distribution and \cref{eq:sde_analytical_solution}. We compute the solution for $CTMGRIT^{loc}$ using the sparse grid level $L=20$ with the initial level $L_0=5$ for the computational domain $[-6, 6]^4$. We run the simulation for $5000$ time steps until time $T_\text{end} = 20$. Here, a local domain size of $S=16$ and a local coarse grid size of $q=12$ per process was used. \cref{fig:linear_oscillator_4d_origin} shows a comparison of the computed solution with the analytical solution. We note again a good agreement between the computed and the analytical solution, with the relative error mostly below five percent over time. The increase in the error compared to the two dimensional case can be attributed to the spatial resolution, which due to the hardware limitations and the choice of temporal parallelization factor $\hat{P}^t_l = 5$ is constrained to be small enough to fit on at most $P^x_l \leq 576$ spatial processes. Hence, the initial level $L_0=5$ has been chosen one level coarser than in the two dimensional case. We expect, in accordance to the theory, that the four dimensional case with an initial level $L_0=6$ will exhibit the same errors as the two dimensional case, due to the fact that increasing the level by one halves the grid spacing, which quarters the error of the computed solution.

\section{Concluding remarks}
\label{sec:concluding_remarks}

In this article, we combined \texttt{XBraid}'s \texttt{MGRIT} parallel-in-time integrator \cite{xbraid-package} with our sparse grid combination approach \cite{Griebel1992} and our space-filling curve based domain decomposition linear solver \cite{Griebel2023}. The resulting overall solver for parabolic problems thus exhibits multiple levels of parallelism: a large number of independent (i.e.\ embarrassingly parallel) subproblems of the combination method in the spatial dimensions, the parallel-in-time parallelism in the time dimension for each independent parabolic subproblem, and finally a domain decomposition based parallel linear solver for the remaining elliptic part of the subproblems. Therefore, the presented approach can in principle utilize a huge number of cores simultaneously where only a moderate number of cores participate in the parallel solution of the considered parabolic problem using the MGRIT time integrator on the respective combination method subproblem (which employs an anisotropic grid in space). Thus, the proposed approach can be employed in a load balanced way on very large supercomputers and on moderate-sized clusters for the same problem, i.e.\ the global requirements on hardware resources are rather small while on much larger hardware, if available, resources can still be efficiently utilized with near perfect speedup. The results of our numerical experiments clearly show these good speedup and scale-up properties of the presented approach in up to six spatial and one temporal dimension. A further improvement of the overall performance of the proposed scheme is conceivable via a joint space-and-time coarsening, i.e. a time-adaptive spatial discretization, and is subject of future research.

\section*{Acknowledgments}

Michael Griebel was supported by the \emph{Hausdorff Center for Mathematics} (HCM) in Bonn, funded by the Deutsche Forschungsgemeinschaft (DFG, German Research Foundation) under Germany's Excellence Strategy -- EXC-2047/1 -- 390685813 of the Deutsche Forschungsgemeinschaft.

\FloatBarrier
\bibliography{bibliography}{}
\bibliographystyle{siamplain}
\end{document}